\documentclass[a4paper,11pt]{amsart}
\usepackage[english]{babel} 
\usepackage[leqno]{amsmath}
\usepackage{amssymb,amsthm}
\usepackage{amscd}
\usepackage{enumerate}

\newtheorem{thm}{Theorem}[section]
\newtheorem{lemma}[thm]{Lemma}
\newtheorem{lemmadef}[thm]{Lemma - Definition}
\newtheorem{prop}[thm]{Proposition}

\newtheorem{cor}[thm]{Corollary}
\newtheorem{exstat}[thm]{Example}

\theoremstyle{definition}
\newtheorem{defi}[thm]{Definition}
\newtheorem{nota}[thm]{}

\theoremstyle{remark}
\newtheorem{remark}[thm]{Remark}

\newtheorem{example}[thm]{Example}

\newcommand{\la}{\longrightarrow}
\newcommand{\ha}{\hookrightarrow}

\newcommand{\ov}{\overline}

\def\X{\mathcal X}
\def\Y{\mathcal Y}
\def\ZZ{\mathcal Z}
\def\L{\mathcal L}

\def\O{\mathcal O}

\def\NN{\operatorname{N}}

\def\dcg{\Delta _X}

\def\int{M_X}

\def\md{\underline{d}}

\newcommand{\pr}[1]{\mathbb{P}^{#1}}

\newcommand{\A}{\mathbb{A}}
\newcommand{\Z}{\mathbb{Z}}

\newcommand{\Div}{\operatorname{Div}}

\newcommand{\Spec}{\operatorname{Spec}}

\newcommand{\mdeg}{\operatorname{{\underline{de}g}}}

\newcommand{\mgbar}{\ov{M}_g}

\def\pdgbar{\overline{P}_{d,\,g}}

\def\ad{\alpha^d_f}

\def\au{\alpha^1_f}

\def\aub{\overline{\alpha^1_f}}
\def\aXub{\overline{\alpha^1_X}}
\def\sb{{\tilde{B}^d_X}}
\newcommand{\LN}{\ell_r}
\newcommand{\uni}{E}

\newcommand{\dX}{\dot{\X}}
\newcommand{\dpi}{\dot{\pi}}
\newcommand{\dfib}{\dot{X}}

\newcommand{\QQ}{\mathcal{Q}}

\newcommand{\nq}{q_f}

\def\pX{P^d_X}

\def\pu{P^1_X}

\def\pXb{\overline{P^d_X}}
\def\pub{\overline{P^1_X}}
\def\pf{P^d_f}
\def\pfb{\overline{\pf}}
\def\nf{N^d_f}
\def\pfu{P^1_f}
\def\pfub{\overline{pfu}}

\def\pfb{\overline{P^d_f}}
\def\pfub{\overline{P^1_f}}

\def\bal{\hat{\mu}}

\newcommand{\Pic}{\operatorname{Pic}}

\newcommand{\picf}[1]{\Pic_f^{#1}}

\newcommand{\picX}[1]{\Pic^{#1}X}

\def\pdbst{\overline{\mathcal{P}}_{d,g}}
\def\pdst{{\mathcal{P}_{d,g}}}
\def\mgbst{\overline{\mathcal{M}}_g}

\newcommand{\tw}{\operatorname{Tw}}
\newcommand{\twf}{\tw_f}

\newcommand{\twX}{\operatorname{Tw}\negthickspace_fX}

\newcommand{\lad}{{\hat{X}}}
\newcommand{\pren}{\overline{r}}

\newcommand{\gen}{\X_K}

\begin{document}

\title[On Abel maps of stable curves]
{On Abel maps of stable curves}
\author[Lucia Caporaso and Eduardo Esteves]
{Lucia Caporaso and Eduardo Esteves}\begin{abstract}
We construct Abel maps for a stable curve $X$. Namely,
for each one-parameter deformation of $X$ to a smooth curve,
having regular total space, 
 and each integer 
$d\geq 1$, we construct by specialization a map 
$\alpha^d_X:\dot{X}^d\to\pXb$, where $\dot{X}\subseteq X$ 
is the smooth locus, and $\pXb$ is 
the moduli scheme for balanced 
line bundles on semistable curves over $X$. For $d=1$, 
we show that $\alpha^1_X$ extends to a map 
$\aXub:X\to\pub$, and does not 
depend on the choice of the deformation. 
Finally,
we give a precise description of when $\aXub$  
is injective.

\end{abstract}

\maketitle

\

{\small \tableofcontents}

\section{Introduction}
 
Let $C$ be a smooth projective curve, and $\Pic^dC$
its degree-$d$ Picard variety 
  parametrizing
 line bundles of degree $d$ on $C$.
For each $d>0$ there exists a remarkable morphism, 
often called the $d$-th {\it Abel map}:
$$
\begin{array}{ccc}
C^d&\la &\Pic^dC\\
(p_1,\ldots,p_d)&\mapsto&\O_C(\textstyle\sum p_i).
\end{array}
$$

Such a map 
has been extensively studied and used in the 
literature. For example
if $d=1$, after  a base point on $C$ is chosen, 
it gives  the Abel--Jacobi embedding 
$C\ha \Pic^1 C\cong \Pic^0 C$  (unless $C\cong\pr{1}$). 
For an interesting historic survey see \cite{K04} or \cite{K05}.

What about Abel maps 
for
singular curves? 
Not much is known for reducible curves,
whereas the case of integral curves is better understood,
as we will explain 
further down.
In the present paper we  construct Abel maps for 
 stable curves 
(including reducible ones, of course).

As we see it, Abel maps should
satisfy the following natural properties.
First, they  should have a 
geometric meaning. More explicitly, recall that for a 
smooth curve $C$ the $d$-th Abel map is 
the moduli map of a natural line bundle on 
$C^d\times C$; see \ref{abelmap}. 
We want a similar property to 
hold for singular curves as well.

Second, 
Abel maps should vary 
continuously in families. In particular, given a 
one-parameter family of smooth curves specializing to a 
singular one, we expect the $d$-th 
Abel maps of the smooth fibers
to specialize to the $d$-th Abel map of the singular fiber.

Both requirements turn out to be nontrivial. In order to 
address the second one we view stable curves as 
limits of smooth ones.
More precisely, let $X$ be a stable curve. Let
$f:\X\to B$ be a family of 
curves over a local one-dimensional 
regular base $B$, with regular total space $\X$, 
smooth generic fiber and $X$ as closed fiber.
 We observe that there exists a 
canonical way to partially extend the $d$-th Abel
map of the generic fiber of $f$, by using 
the N\'eron model $N^d_f$ of the degree-$d$ Picard scheme 
of that fiber. The N\'eron mapping property yields a 
close relative of the $d$-th Abel map of $X$, 
defined on the nonsingular locus 
${\dot X}^d\subseteq X^d$, 
which 
we call the $d$-th {\it Abel--N\'eron map} of $X$; see 
\ref{key}. The target of 
this
map is the closed fiber of 
$N^d_f$, rather than the Picard scheme.

The great advantage of Abel--N\'eron maps is their 
naturality, obtained directly from
the universal property of N\'eron models. However, they have
 two major drawbacks. First, they
do not have any a-priori modular interpretation. Second, 
they are not defined on 
the whole 
$X^d$.

To attack these problems we consider the modular 
compactified Picard variety introduced in 
\cite{caporaso} and further studied in \cite{cner}. 
If $X$ is suitably general,
more precisely  ``$d$-general" 
(Definition \ref{dgen}),
 there exists  a natural proper 
$B$-scheme $\pfb$   parametrizing equivalence 
classes of ``balanced" 
line bundles on semistable curves; see \ref{pdg}. 
These are line bundles whose multidegree satisfies certain 
inequalities;
see 
Definition \ref{baldef}.
It is shown 
in \cite{cner} that 
$\pfb$ contains $N^d_f$ as a dense open subscheme, hence
it 
gives
a geometrically meaningful  description and  completion of $N^d_f$. 

So, assume for now
that $X$ is $d$-general.
 Let $\pXb$ be the  closed fiber of $\pfb$; it 
does not depend on 
the choice of 
$f$. Viewing $N^d_f$ inside $\pfb$, it is thus possible to 
exhibit the $d$-th Abel--N\'eron 
map as a map $\alpha^d_X:{\dot X}^d\to\pXb$ (see Proposition \ref{keyrig}) and thus 
give it a 
geometric interpretation.
Moreover, since $\pXb$ is 
complete, we 
are able to study whether $\alpha^d_X$ extends.

Although $\alpha^d_X$ is modular, an explicit 
description for it 
is hard to exhibit in full generality; we do that only for curves with two components (see
Proposition~\ref{vine}).  Without such a description it is difficult to 
study whether and how $\alpha^d_X$ extends.

The case $d=1$ turns out to be easier. By means of 
Theorem~\ref{AN1} we give an explicit description 
of the line bundle defining $\alpha^1_X$. 
Using 
this description,
we show in Corollary \ref{AX}
that $\alpha^1_X$ does not depend on 
the choice of 
$f$, a remarkable 
property not to be expected  in general for 
$d>1$; see Remark~\ref{vinermk}.

Refining our
modular description of $\alpha^1_X$, we  construct a completion for it as a regular map
$\aXub:X\to\pub$ 
(Theorem \ref{A1c}).
Finally, we prove
that $\aXub$ 
is as close as it can be to an injection
(see Proposition~\ref{A1cinj} for the precise statement). 

Now suppose that $X$ is not 
1-general,
so that $X$ lies in a proper closed subset of $\mgbar$
for $g$ even; see Proposition~\ref{1gen}.
Then
$\pfub$ fails to contain $N^1_f$; nevertheless our existence results do
extend, suitably modified (see 
\ref{nogen}),
whereas uniqueness and injectivity results (like Proposition~\ref{A1cinj}) 
may fail. 
In
this case
the 
setup
is significantly more complicated
for standard technical reasons (presence of 
non-GIT-stable points, or of nonfine moduli spaces.) This is why we
chose to  first work  under the assumption of $1$-generality,
and to later indicate, in 
\ref{nogen} and \ref{Q1fix},
how to modify  proofs and  statements to include  the
special case.

\

Abel maps 
were  constructed
for all integral  curves
in \cite{AK}, and further studied in \cite{EGK1},
\cite{EGK2} and \cite{EK}.
In \cite{AK}, it is shown that the
first Abel map of an integral singular curve is an embedding into its
compactified
Picard scheme.

 Constructing Abel maps for reducible curves 
presents further difficulties,
due to the lack of natural, separated target spaces.
The use of N\'eron models 
as target spaces
is not new in the literature: in \cite{edix}
Abel--Jacobi maps for nodal curves were 
studied by means of the N\'eron mapping property,
similarly to what we do here with our 
Abel--N\'eron maps. 
However, N\'eron models are seldom proper
 and 
thus we cannot expect Abel maps to N\'eron models to be 
defined everywhere. In this framework, our
contribution 
is that
of bringing compactified Picard schemes into the picture.
This enables us to compactify N\'eron models 
and hence to obtain a target space 
into which complete Abel maps could be defined. In fact, 
we prove that $\alpha^1_X$ extends;  it is still an open 
problem how to extend $\alpha^d_X$  for $d>1$.

It is our intention to 
construct the 
first Abel map of any Gorenstein singular curve, 
using the 
modular compactified Picard schemes introduced in 
\cite{Est}.  
For such curves 
the connection with N\'eron models is not available.
Nonetheless, prompted by
the results of 
the present
paper, we found 
 significant evidence for the existence of  a  modular interpretation
similar to the present  case.

Our paper is organized as follows. Section~\ref{nerpic} 
is devoted to preliminaries of various types.
In Section~\ref{balanced} we describe degree-$d$ 
Abel--N\'eron maps 
to the compactified Picard scheme.
In Section~\ref{open} we establish the modular description of the Abel--N\'eron map in 
degree $1$,
and show that it is independent of the choice of the deformation.
Finally, in Section~\ref{close}  we construct the completed degree-$1$ Abel map,
give it a modular description and study when it is injective. Finally, from
\ref{nogen} to the end of the paper, we explain how to handle the special case of non-$1$-general curves.

\section{N\'eron models of Picard schemes}
\label{nerpic}
\begin{nota}
\label{not}
{\it Setup.}
We work over a fixed algebraically closed field $k$. 
All schemes are assumed locally of 
finite type over $k$, unless stated otherwise

For us, a {\it curve} is a reduced and 
connected projective scheme of 
 dimension 1.
Mostly, we will deal 
with {\it nodal curves}, that is, curves whose only 
singularities are nodes.

A {\it regular pencil} ({\it of curves}) 
is a flat projective morphism $f:\X\to B$ 
between connected, regular,
schemes such that 
$\dim B=1$, every geometric fiber of $f$ is a curve, and 
$f$ is smooth over a dense open subscheme of $B$. 

We call a regular pencil $f:\X\to B$ {\it local} if 
$B=\Spec R$, where $R$ is a discrete valuation ring 
(having $k$ as residue field). If $X$ is 
the closed 
fiber, we will also say that $f$ is a {\it regular 
smoothing of $X$}. 

For each regular pencil $f:\X\to B$ we shall let 
$K:=k(B)$, the field of rational functions of $B$, 
and denote by $\gen$ the generic fiber of $f$. Notice 
that $\gen$ is a smooth curve over $K$.

Given any morphism $f:\X \to B$, and any integer 
$d\geq 1$, let 
$f_d:\X^d_B\to B$ denote 
the $d$-th fibered power of $\X$ over $B$. 
If $f$ is a regular pencil, 
we denote the open subset of $\X^d_B$ where $f_d$ 
is smooth by $\dX^d_B$; so
$$
\dX^d_B:=\X^d_B\smallsetminus\text{Sing}(f_d)
$$
If $f$ is a local regular pencil, let 
$\dfib ^d$ denote the closed fiber of 
$\dX^d_B\to B$; so
$$
\dfib^d = \{(p_1,\ldots,p_d): p_i\in X\smallsetminus 
X_{\text{sing}}\}.
$$

Given any morphism $f:\X \to B$ and any 
$B$-scheme $T$, the base change of $f$ to $T$ is
denoted $f_T:\X_T\to T$. 
\end{nota}

\begin{nota}
{\it The relative Picard scheme.}\label{modmapL} 
Let $f:\X\to B$ be a regular 
pencil, and $d$ an integer. 
The closed fibers of $f$ are 
geometric, by our general assumption, and the general 
fiber is smooth. Thus the irreducible components of the 
fibers of $f$ are geometrically irreducible. By 
a theorem of Mumford's, \cite{BLR}, Thm. 2, p.~210, 
the (relative) Picard scheme $\Pic_f$ of $f$ exists, 
and is locally of finite type over $S$.
Furthermore, $\Pic_f$ is 
formally smooth over $B$ by \cite{BLR}, Prop. 2, p.~232, 
whence smooth over $B$ by \cite{BLR}, Prop. 6, p. 37. 

Let $\picf{d}$ be the degree-$d$ 
Picard scheme of $f$, the open subscheme of $\Pic_f$ 
  parametrizing line bundles of relative degree 
$d$. Given any $B$-scheme $T$ and any 
line bundle $\L$ on $\X_T$ of $f_T$-relative degree $d$, 
there is a moduli map associated to $\L$,
\begin{equation}
\begin{array}{lccr}
\mu_\L: & T&\la &\picf{d}\\
&t &\mapsto &L_t,
\end{array}
\end{equation}
where $L_t\in \Pic^d X_t$ is the restriction of 
$\L$ to the fiber $X_t:=f_T^{-1}(t)$. 
The map $\mu_\L$ determines $\L$ up to 
tensoring with pullbacks of line bundles
from $T$. Notice that to a map $T\to \picf{d}$ there 
does not necessarily correspond a line bundle on 
$\X_T$, though the line bundle will exist, for instance, 
if $f$ admits a section; see \cite{BLR}, Prop. 4, p. 204.
\end{nota}

\begin{nota}
\label{Npre}
{\it N\'eron models of Picard schemes.} Let 
$f:\X\to B$ be a regular pencil, and $d$ an integer. 
Recall that a 
basic characteristic (and a drawback for various 
applications) of the Picard scheme $\Pic^d_f$ is that it 
is not separated over $B$, if $f$ has 
reducible special fibers. One way to fix this is to
introduce the N\'eron model:
     $$N^d_f:=\NN(\Pic^d\gen).$$ 

The N\'eron model is a smooth, separated 
(possibly not proper) scheme of finite type over $B$ 
with generic fiber equal to $\Pic^d\gen$, 
which satisfies a fundamental mapping 
property that
uniquely determines 
it. Namely, for every smooth $B$-scheme $Z$ each map 
$Z_K\to\Pic^d\gen$ extends uniquely to a map 
$Z\to N^d_f$; see \cite{BLR}, Def.~1, p. 12. 

The existence of $N^d_f$ for any regular pencil $f$ is 
likely well known.
 Since this 
result
is fundamental for our work, but we could not find the
precise statement to refer to, 
we sketch 
a proof of it using results in \cite{BLR}.
First, assume that $f$ is local, that is, 
$B$ is the spectrum of a discrete valuation ring $R$. 
Then there is a N\'eron model of $\Pic^d\gen$ over $B$, 
which is equal to $\Pic^d_f$ if $f$ is smooth. Indeed, 
since $\Pic^d\gen$ is a $\Pic^0\gen$-torsor, 
by descent theory we may assume that 
$R$ is a strictly Henselian ring; see 
\cite{BLR}, Cor. 3, p. 158. In this case, 
$f$ admits a section through its smooth locus, by 
\cite{BLR}, Prop. 5, p.~47. This section can be used 
to produce a $B$-isomorphism $\Pic^d_f\to\Pic^0_f$. 
We may thus assume that $d=0$. In this case, 
there is a N\'eron model of $\Pic^0\gen$ over $B$ because 
$\Pic^0\gen$ is an Abelian variety over $K$; see 
\cite{BLR}, Cor. 2, p. 16. Furthermore, if $f$ is smooth, 
$\Pic^0_f$ is an Abelian $B$-scheme, whence is the 
N\'eron model of $\Pic^0\gen$ over $B$ by \cite{BLR}, 
Prop. 8, p. 15.

Now, consider the general case. Let $U\subseteq B$ 
be the largest open subscheme over which $f$ is smooth, 
and set $h:=f|_{f^{-1}(U)}$. As we saw above, 
$\Pic^d_h$ restricts to the N\'eron model of 
$\Pic^d\gen$ locally around each point of $U$. Then 
$\Pic^d_h$ is the N\'eron model of $\Pic^d\gen$ over $U$ 
by \cite{BLR}, Prop. 4, p. 13. Finally, the local 
existence of N\'eron models of $\Pic^d\gen$ around each 
point of $B$, and the existence of the N\'eron model 
over the dense open subscheme $U\subseteq B$ imply the 
(global) existence of the N\'eron model over $B$, by 
\cite{BLR}, Prop. 1, p. 18. 

Since $\picf{d}$ is smooth over $B$, a 
first consequence of the mapping property of the 
N\'eron model $N^d_f$ is the 
existence of a canonical $B$-morphism
\begin{equation}
\label{quotmap}
\nq:\picf{d}\la N^d_f
\end{equation}
which is the identity on the generic fiber.

Assume now that the geometric fibers of $f$ are 
nodal. Let $X$ be a closed fiber of $f$. In the 
description of the N\'eron model, and 
also in our paper, the following set 
$\twf X$ of (isomorphism classes of) distinguished 
line bundles plays an important role:
$$
\twf X:=
\frac{\{\O_\X(D)|_X:D\in\Div\X\text{ with }
\text{Supp }D\subset X\}}{\cong}
\subset\Pic^0X.
$$
The divisors $D$ appearing above are simply sums with 
integer coefficients of the components of $X$, which are 
Cartier divisors of $\X$ 
because $\X$ is regular. Line bundles 
in $\twf X$  are called {\it twisters}. 
Here is a useful observation:
\begin{equation}
\label{degtwi}
\forall\, T,T'\in \twf X,\hskip.3in T=T'\Leftrightarrow
\mdeg T=\mdeg T'.
\end{equation}
Since twisters are specializations of the 
trivial line bundle of the generic fiber, $\O_{\gen}$, 
all of them must be identified in any
separated quotient of $\Pic^0_f$. In particular, 
$\nq(T)=\nq(\O_X)$ for each $T\in\twf X$. 

We shall now identify multidegrees that differ by 
multidegrees of twisters. Let $\gamma$ be the 
number of irreducible components of $X$, and set 
$$
\Lambda_X:=\{\mdeg T:T\in \twf X\}\subseteq\Z^\gamma.
$$
Define now an equivalence relation ``$\equiv$" on 
multidegrees by setting
$$
\md \equiv \md' \Leftrightarrow \md -\md'\in \Lambda_X.
$$
The set of multidegree classes $\md+\Lambda_X$ 
with fixed total 
degree $|\md|$ equal to $d$ is denoted by $\dcg^d$. Thus
\begin{equation}
\label{dcg}
\dcg^d:=\frac{\{\md \in \Z^{\gamma}: |\md |=d\}}{\equiv}.
\end{equation}

It is well known that $\dcg^0$ is a finite group, a 
purely combinatorial invariant of $X$ 
often called the  degree class group of $X$,
or the group of connected components of $N^d_f$. 
In addition, for each $d$ there is a (nonunique) 
bijection $\dcg^0\to\dcg^d$, obtained by summing with any 
multidegree $\md$ with $|\md|=d$.

For each $\delta\in\dcg^d$, let $\md$ be any multidegree 
representing $\delta$, and set
\begin{equation}
\label{picdel}
\picf{\delta}:=\picf{\md}=
\{\L\in \Pic^d_f: \mdeg\L|_X = \md\}\subset\Pic^d_f.
\end{equation}
This definition is easily seen not to  depend on the choice of the 
representative $\md$ (cf. \cite{cner}, 3.9). 

Assume now that $f$ is a regular smoothing of $X$.
At this point we are able to describe the N\'eron model of
$\Pic^d\gen$:
\begin{equation}
\label{nerdesc}
\NN^d_f\cong\frac{\coprod _{\delta \in \dcg^d}
\picf{\delta}}{\sim_K},
\end{equation}
where $\sim_K$" denotes the gluing along the generic 
fiber, equal to $\Pic^d\gen$; 
see \cite{cner}, Lemma 3.10.

Let $N^d_X$ denote the closed fiber of $N^d_f$. Observe 
that $N_X^d$ is a disjoint union of finitely many copies 
of the generalized Jacobian of $X$:
picking a representative $\md^{\delta}$ for each class 
$\delta\in\dcg^d$, we have
$$
N_X^d\cong\coprod_{\delta\in\dcg^d}\picX{\md^\delta}.
$$
Although the above isomorphism is not canonical, 
we see that the scheme structure of $N_X^d$ does
not depend on $f$. The closed points of $N^d_X$ are 
in 1-1 correspondence with the degree-$d$ line 
bundles on $X$ modulo twisters. In  particular, 
for $d=0$, we have  
$
\nq^{-1}(\nq (\O_X))=\twX
$.
\end{nota}

\begin{nota}
\label{nosep} {\it N\'eron maps.}
Let $f:\X\to B$ be a regular pencil. 
Let $T$ be a $B$-scheme, and $\L$ a line bundle on 
$\X_T$ of relative degree $d$ over $T$. Let 
$\mu_\L:T\to\picf{d}$ be the moduli map of $\L$, defined 
in \ref{modmapL}.
Consider the composition:
$$
\overline{\mu_{\L}}:T\stackrel{\mu_{\L}}{\la}\picf{d}
\stackrel{\nq}{\la} N^d_f.
$$
We call $\overline{\mu_{\L}}$ the {\it N\'eron map of} 
$\L$. Notice that $\L$ is certainly not 
determined by its  N\'eron map,
not even modulo pullbacks of line bundles on $T$. 
In fact, if $D\subset\X$ is a Cartier divisor 
entirely supported on a closed fiber of $f$, then 
$\L\otimes \O_{\X_T}(D_T)$ has the same N\'eron map as 
$\L$, because $\nf\to B$ is separated. 
\end{nota}

\begin{nota}{\it Abel--N\'eron maps.}
\label{abelmap} 
Let us recall the precise definition of the Abel map of a smooth curve, using the same
set up of   \cite{GIT}, Section 6 p.118, 119. 
 Let
$
h:{\mathcal C}\to S
$
be a smooth curve over a scheme $S$,
so that $h$ is a smooth 
morphism  whose  fibers 
are curves. 
For each integer $d\geq 1$, let 
${\mathcal C}^d_S$ be the $d$-th fibered power of
$\mathcal C$ over $S$. There 
is a canonical $S$-morphism
\begin{equation}
\begin{array}{lccr}
& {\mathcal C}^d_S&\la &\Pic_h^d \hskip.5in  \\
& C_s^d \ni (p_1,\ldots,p_d) &\mapsto 
&\O_{C_s}(p_1+\cdots+p_d),
\end{array}
\end{equation}
defined over each $s\in S$ by taking a $d$-tuple of points 
of the fiber $C_s$ to the line bundle associated to their 
sum, which we shall call the {\it $d$-th Abel map} 
of $h$. 
Recall 
that the above map is the moduli map of
a natural line bundle on
${\mathcal C}^d_S \times _S {\mathcal C}$, namely the one 
associated to the Cartier divisor 
$\sum _1^d S_i$, where each $S_i$ is the image of the
$i$-th natural section $\sigma _i$ of the first projection
${\mathcal C}^d_S\times_S{\mathcal C}\to
{\mathcal C}^d_S$, given by
\begin{equation}
\label{natsec}
\sigma_i(p_1,\ldots,p_d)=((p_1,\ldots,p_d),p_i).
\end{equation}

We may apply this construction to 
a regular pencil $f:\X\to B$. First of all, since 
$\X_K$ is smooth over $K$, we may consider 
the $d$-th Abel map of $\X_K$:
\begin{equation}
\label{abgen}
\alpha^d_K : \gen^d \la \Pic^d\gen.
\end{equation}
The above map extends to a map $\dX^d_B\to\Pic_f^d$. 
Indeed, the extension is the moduli map of the line bundle 
associated to the Cartier divisor $\uni^d$ of 
$\dX^d_B\times_B\X$, where $\uni^d$ is the sum of the 
images of the $d$ natural sections 
$\sigma_1,\ldots,\sigma_d$ given by (\ref{natsec}) 
of the first projection 
$\dX^d_B\times_B\X\to\dX^d_B$. Composing with 
the map $q_f$ of \eqref{quotmap}, we obtain the 
N\'eron map of $\O_{\dX^d_B\times _B\X}(\uni^d)$ 
(cf. \ref{nosep}), which is also an extension of 
$\alpha^d_K$. 
\end{nota}

The first simple but crucial observation is the following 
(well-known) fact:

\begin{lemmadef}
\label{key} Let $f:\X\to B$ be a regular pencil. 
For each integer $d\geq 1$ there exists a 
unique morphism, which we call the 
{\em $d$-th Abel--N\'eron map of $f$},
$$
\NN(\alpha ^d_K): \dX^d_B \la \NN(\Pic^d\gen)=N^d_f
$$
whose restriction to the generic fiber is 
$\alpha^d_K$. The map $\NN(\alpha ^d_K)$ 
is the N\'eron map of $\O_{\dX^d_B\times _B\X}(\uni^d+D)$ 
for every Cartier divisor $D$ of $\dX^d_B\times_B\X$ 
supported on any finite number of closed fibers 
of $\dX^d_B\times_B\X\to B$.
\end{lemmadef}

\begin{proof} The existence of an extension 
of $\alpha^d_K$ 
is a straightforward consequence of the N\'eron 
mapping property, \cite{BLR}, Def. 1, p. 12:
since $\dX^d_B $ is smooth over $B$ and has 
generic fiber $\gen^d$, the Abel map 
$\alpha ^d_K$ admits a 
unique extension 
$\NN(\alpha ^d_K):\dX^d_B\to\NN(\Pic^d\gen)$.

Since also the N\'eron map of 
$\O_{\dX ^d_B\times _B\X}(\uni^d+D)$, for any 
$D$ as described in the lemma, 
extends $\alpha^d_K$, the last statement 
follows from the fact that $N^d_f$ is separated over $B$. 
\end{proof}

\section{Abel maps to balanced Picard schemes}
\label{balanced}
We want to give a modular interpretation of 
the Abel--N\'eron maps and, at the same time, study the 
problem of completing them. To do this we shall use 
some results of \cite{cner}, where  N\'eron models 
are glued together over the moduli space of stable
curves and are thereby endowed of a 
geometrically meaningful completion.  
Our moduli problem is centered around 
Definition~\ref{baldef} below. First, we recall 
a few concepts.

\begin{nota}
\label{notcurve}
Let $X$ be a nodal curve of arithmetic 
genus $g\geq 2$. Denote by 
$\omega_X$ its canonical, or dualizing bundle. 
For each proper subcurve $Z\subsetneq X$, 
always assumed to be complete, let 
$Z':=\overline{X\smallsetminus Z}$ and 
$k_Z:=\#Z\cap Z'$. Also, let $w_Z:=\deg_Z\omega_X$. 
If $Z$ is connected, denote by $g_Z$ 
its arithmetic genus, and recall that 
\begin{equation}\label{adjw1}
w_Z=2g_Z-2+k_Z,
\end{equation}
a well-known identity that can be proved using 
adjunction.

We call $X$ {\it semistable} 
(resp. {\it stable}) if $k_Z\geq 2$ (resp. $k_Z\geq 3$) 
for each smooth rational component $Z$ of $X$. 
Those $Z$ for which $k_Z=2$ are called {\it exceptional}.
A semistable curve is called {\it quasistable} if 
two exceptional components never meet each other. If 
$X$ is semistable, it follows from \eqref{adjw1} 
that $w_Z\geq 0$ for each subcurve 
$Z\subseteq X$, with equality if and only if $Z$ 
is a union of exceptional components.

A {\it family of semistable} (resp. {\it stable}, resp. 
{\it quasistable}) {\it curves} is a flat, projective map 
$f:\X\to B$ whose geometric 
fibers are semistable (resp. stable, 
resp. quasistable) curves. A line bundle of degree $d$ 
on such a family $f:\X\to B$ is a line bundle on $\X$ 
whose restriction to each fiber has degree~$d$. 
\end{nota}

\begin{defi}
\label{baldef} Let $X$ be a semistable curve 
of arithmetic genus $g\geq 2$, and let $L\in \Pic^dX$.
\begin{enumerate}[(i)]
\item
\label{sem}
We say that $L$ and its multidegree $\mdeg L$ are
{\it semibalanced} 
if for each connected proper subcurve 
$Z\subsetneq X$ the 
{\it Basic Inequality} below holds:
\begin{equation}
\label{BIeq}
m_Z(d)\leq \deg_ZL\leq M_Z(d),
\end{equation}
where
$$
M_Z(d):=\frac{dw_Z}{2g-2}+\frac{k_Z}{2} 
$$
and
\begin{displaymath}
m_Z(d):=\left\{ \begin{array}{ll}
M_Z(d)-{k_Z} 
&\text{ if } Z \text{ is not an exceptional component,}\\ 
0 &\text{ if } Z \text{ is  an exceptional component.}\\
\end{array}\right.
\end{displaymath}
\item
\label{bal}
We call $L$ and $\mdeg L$ {\it balanced} if they 
are semibalanced and if for each exceptional 
component $E\subset X$ we have
$$
\deg_E L=1.
$$
\item
We call $L$ and $\mdeg L$ {\it stably balanced} 
if they are balanced and if for each connected 
proper subcurve $Z\subsetneq X$ 
such that $\deg _ZL=m_Z(d)$ the complement 
$Z'$ is a union of exceptional 
components.
\item
\label{balset} We denote by $B^d_X$ the set of 
balanced multidegrees on $X$, 
and by $\sb$ the subset of stably balanced 
ones.\footnote{In \cite{cner}
$B^d_X$ denotes the semibalanced multidegrees.}
\item
A line bundle (of degree $d$) 
on a family of semistable curves 
is called semibalanced (balanced or stably balanced)
if its restriction to each geometric 
fiber of the family is.
\end{enumerate}
\end{defi}

\begin{remark}
\label{balrm} 
Definition~\ref{baldef} is 4.6 of \cite{cner}, 
which originates from 
\cite{caporaso}, Section 3.1 and from \cite{CCC}, 5.1.1. There the Basic Inequality is 
$$
\Big|\deg_ZL-\frac{dw_Z}{2g-2}\Big|\leq\frac{k_Z}{2},
$$
and the extra condition for when $Z$ is an exceptional 
component is imposed separately. For convenience, we 
presented a set of inequalities, \eqref{BIeq}, 
that includes the exceptional cases. Abusing the 
terminology, we still call \eqref{BIeq} the 
Basic Inequality.

We mention some simple but useful consequences of the 
definiton.
\begin{enumerate}[(A)]
\item
\label{balrm1}
If $X$ is stable, then (\ref{sem}) and (\ref{bal}) 
coincide, i.e. a semibalanced line bundle is always 
balanced. 
\item
\label{balrm2}
There can only exist a balanced line bundle on 
a semistable curve $X$ if $X$ is quasistable.
Indeed, let $Z\subset X$ be a 
connected chain of exceptional components. If $L$ is a 
semibalanced line bundle on $X$, then $L$ has degree~0 
on every component of $Z$ but possibly one, where 
$L$ may have degree 1. 
However, if $L$ is balanced, $L$ cannot have 
degree 0 on any component of $Z$. Then 
$Z\cong \pr{1}$ and $\deg_ZL=1$. 
\item
\label{balrm3} To check whether $L$ is semibalanced it is 
enough to check whether $\deg_ZL\geq m_Z(d)$ for 
each connected proper subcurve $Z\subsetneq X$. 
Indeed, let $Z$ be such a subcurve, and 
let $Y_1,\dots,Y_n$ denote the 
connected components of $Z'$. By hypothesis,
$$
\deg_{Y_i}L\geq m_{Y_i}(d)\geq\frac{dw_{Y_i}}{2g-2}-
\frac{k_{Y_i}}{2}
$$
for each $i$. Since 
$
k_{Y_1}+\cdots+k_{Y_n}=k_{Z'}=k_Z,
$
summing up the above inequalities we get
$$
\deg_{Z'}L\geq \frac{dw_{Z'}}{2g-2}-
\frac{k_Z}{2}.
$$
Now, since $\deg_ZL+\deg_{Z'}L=d$ and 
$w_Z+w_{Z'}=2g-2$, we get
$$
\deg_Z L=d-\deg_{Z'}L\leq d-\frac{dw_{Z'}}{2g-2}+
\frac{k_Z}{2}=\frac{dw_Z}{2g-2}+\frac{k_Z}{2}.
$$
\end{enumerate}
\end{remark}

\begin{nota}
\label{sbrep} In \cite{cner}, Lemma 4.4, it is proved that 
each multidegree class has a semibalanced representative.
More precisely, fix an integer $d$, and let $X$ be a 
stable curve. Recall the notation in \eqref{dcg} and 
\ref{baldef} (\ref{balset}). Then 
Lemma 4.4 of \cite{cner} implies that the natural map 
below is surjective (square brackets denoting classes):
\begin{equation}
\label{mapclass}
\begin{array}{lccr}
[\ldots]: &B^d_X&\la &\dcg^d \\
&\md &\mapsto &[\md]\  
\end{array}
\end{equation}

We shall say that $X$ is ``$d$-general'' if 
the map (\ref{mapclass}) is bijective; see 
Definition \ref{dgen} below
\end{nota}

\begin{nota}
\label{C94}
The moduli problem for balanced line bundles was 
introduced and studied in \cite{caporaso}
to compactify the universal Picard scheme over $M_g$.
That compactification was constructed as a GIT-quotient. 
We do not need to recall the details of the construction 
here, only a few facts.
There are morphisms
$$
H_d \stackrel{\pi_d}{\la}H_d/G\stackrel{\phi_d}{\la}\mgbar
$$
where $H_d$ is an open subscheme of a suitable Hilbert 
scheme, acted upon by an algebraic group $G$, the map 
$\pi_d$ is a GIT-quotient map, and $\phi_d$ is a 
surjective, projective morphism. The quotient scheme 
$H_d/G$, denoted by $\pdgbar$, 
is integral and projective. The fiber of
$\phi_d$ over a general smooth curve $X$ is exactly 
$\Pic^dX$; see \cite{caporaso}, Thm. 6.1.

Denote by $U\subseteq\pdgbar$ the nonempty open subscheme 
over which $\pi_d$ restricts to 
a geometric GIT-quotient, i.e. all fibers are orbits, and 
all stabilizers are finite and reduced. 

\begin{defi}
\label{dgen} Let $X$ be a stable curve of arithmetic 
genus $g\geq 2$. 
We say that $X$ is {\it $d$-general}
if any of the following equivalent conditions hold:
\begin{enumerate}[(i)]
\item
$\phi_d^{-1}(X)\subset U$.
\item
The class map (\ref{mapclass}) of \ref{sbrep} is bijective.
\item
Every balanced 
line bundle on $X$ is stably balanced, i.e. $\sb=B^d_X$.
\end{enumerate}
\end{defi}
The  equivalence of these three conditions follows 
from \cite{caporaso} Lemma 6.1. 

It is known that 
{\it every stable curve of genus $g$ is $d$-general
if and only if the integers $d$ and $g$ satisfy 
$(d-g+1,2g-2)=1$}; see \cite{caporaso}, Prop 6.2.
\end{nota}

\begin{nota}
\label{equiv}
We need to recall when two semibalanced line bundles 
are defined to be equivalent. Let $X$ be a stable curve, 
and $X_1$ and $X_2$ two semistable curves having $X$ as 
stable model. For each $i=1,2$ let $L_i$ be a 
semibalanced line bundle on $X_i$. 
Let $Y_i$ be the semistable curve 
obtained by contracting all exceptional components 
of $X_i$ where $L_i$ has degree 0. Then there is a 
unique line bundle $M_i$ on $Y_i$ 
whose pullback to $X_i$ is $L_i$. Since $L_i$ is 
semibalanced, $M_i$ is balanced and $Y_i$ is 
quasistable. Let $F_i\subseteq Y_i$ be 
the union of all the exceptional components of 
$Y_i$, and let $\widetilde{Y_i}:=F'_i$, the 
complementary subcurve of $F_i$ in $Y_i$. 
Then {\it $L_1$ and $L_2$ are  equivalent 
if $Y_1=Y_2$ and 
$M_1|_{\widetilde{Y_1}}\cong M_2|_{\widetilde{Y_2}}$}. 
Notice that $M_i$ is equivalent to $L_i$ for $i=1,2$. 

Thus, every equivalence class includes always a 
balanced line bundle $N$ on a quasistable curve $Y$. 
The quasistable curve $Y$ is unique, but $N$ is not. 
What is unique is the restriction of $N$ to $\tilde Y$, 
the complementary subcurve of the union $F$ of all the 
exceptional components of $Y$. The 
quasistable curve $Y$ and $N|_{\tilde Y}$ determine 
the equivalence class. The restriction $N|_F$ 
is also unique, since a balanced line bundle must 
have degree 1 on every exceptional component.
So, our equivalence relation disregards the 
gluing data of the bundles over the points in 
$\tilde Y\cap F$. 

If $\X\to B$ is a family of 
semistable curves, then two semibalanced line bundles 
$\L_1$ and $\L_2$ on $\X\to B$ are called equivalent 
if and only if their restrictions to every geometric 
fiber of $\X\to B$ 
are equivalent in the sense explained above.
\end{nota}

\begin{nota}
\label{pdg} Let $d$ and $g$ be integers, with $g\geq 2$. 
Assume first that $d-g+1$ and $2g-2$ are coprime, 
so that every stable curve of 
arithmetic genus $g$ is $d$-general.
Then the construction summarized in 
\ref{C94} can be improved, by considering stacks. 
More precisely, there
exist two (modular) Deligne--Mumford stacks
$\pdbst$ and $\pdst$, each one equipped with a natural, 
strongly representable morphism
to $\mgbst$. (To tie in with \ref{C94}, notice that 
$\pdbst$ is the quotient stack $[H_d/G]$.)
The following properties hold;
see \cite{cner}, Section~5 for details:

\begin{enumerate}[{\bf (A)}]
\item
\label{fiber}
For each ($d$-general) stable curve $X$, denote by 
$\pX$ and $\pXb$ the fibers 
of $\pdst$ and $\pdbst$ over $X$. 
Since $\pdst$ and $\pdbst$ are strongly representable over 
$\mgbst$, both $\pX$ and $\pXb$ are quasiprojective 
schemes.

The first, $\pX$, is the fine moduli scheme of degree-$d$ 
balanced line bundles on $X$. 
The second, $\pXb$, is the coarse 
moduli scheme of equivalence classes 
of degree-$d$ semibalanced line bundles on 
semistable curves having $X$ as stable model; see 
\ref{equiv}. Actually, $\pXb$ 
is not far from being a fine moduli scheme; see 
(\ref{fine}) below. 

$\pX$ lies naturally inside $\pXb$ as an open and dense subscheme. 

\item
\label{mymap}
Let $f:\X \to B$ be any family of ($d$-general) 
stable curves of genus $g$, and consider the schemes 
$$\pf:=B\times _{\mgbst } \pdst
\  \  \text{ and }
\   \  \pfb:=B\times _{\mgbst } \pdbst.$$
(That these are indeed schemes follows, again, 
from the fact that the maps $\pdst\to \mgbst$ and 
$\pdbst\to \mgbst$
are strongly representable.)
We have a natural inclusion $\pf \subset \pfb$. 
An explicit description of $\pf$, for when $f$ is a 
local regular pencil, is (\ref{pfeq}). As for $\pfb$, the 
following fact holds: to each 
triple $(T, \Y\to T, \L)$ where 
 $T$ is a $B$-scheme, $\Y \to T$ is a family 
of semistable curves 
having $f_T:\X_T\to T$ as stable model,
and $\L$ is a semibalanced line bundle of degree 
$d$ on $\Y\to T$, there corresponds 
a moduli map
$$
\bal_{\L}:T\la \pfb,
$$
taking each geometric point $t$ of $T$ to the equivalence 
class of the restriction of $\L$ to the (geometric) 
fiber of $\Y$ over $t$. We call $\bal_\L$ 
the {\it moduli map} of $\L$.

The image of $\bal_{\L}$ is contained in $\pf$ if and 
only if $\L$ has degree $0$ on every 
exceptional component of every geometric fiber of 
$\Y\to T$.
\item
\label{fine}
The scheme 
$\pf$ is a fine moduli scheme; see \cite{cner}, Cor. 5.14 
and Rmk.~5.15. 
Also, $\pfb$ is not far from being a 
fine moduli scheme. In fact, it 
is endowed with a quasiuniversal
pair $(\ZZ\to\pfb, {\mathcal N})$, where $\ZZ\to\pfb$ 
is a family of
quasistable curves having 
$f_{\pfb}:\X_{\pfb}\to\pfb$ as stable model, 
and ${\mathcal N}$  is
a  balanced line bundle of degree $d$ on $\ZZ\to\pfb$ 
that has a role similar to that 
of a Poincar\'e bundle. Indeed, 
for each triple $(T, \Y\to T, \L)$ 
where $T$ is a $B$-scheme, 
$\Y \to T$ is a family 
of semistable curves with stable model $f_T:\X_T\to T$,
and $\L$ is a semibalanced line bundle 
of degree $d$ on $\Y\to T$, there is a map 
$\Y\to\ZZ$ such that the diagram of maps below 
is commutative,
$$\begin{CD}
\Y @>>> \ZZ\\
@VVV @VVV\\
T @>\bal_{\L}>> \pfb,\\
\end{CD}$$
and such that $\L$ is equivalent to the pullback of
${\mathcal N}$ to $\Y$; see \ref{equiv}. 
(The map $\Y\to\ZZ$ is certainly 
not uniquely determined, which is why we call the pair 
$(\ZZ\to\pfb, {\mathcal N})$ quasiuniversal.) 
Furthermore, if $\Y\to T$ is a family of quasistable 
curves, and $\L$ is balanced, 
then the map $\Y\to\ZZ$ can be chosen such that 
the above diagram is a fibered product diagram.
\end{enumerate}
\end{nota}

\begin{remark}
\label{noMR}
If $(d-g+1, 2g-2)\neq 1$, almost everything in \ref{pdg} 
works over the open subset of
$\mgbar$   parametrizing $d$-general stable curves. 
For a proof, it suffices to argue exactly as for  Theorem 5.9 in \cite{cner},
after replacing 
$\mgbst$ with the substack of $d$-general  curves,
and the stacks $\pdst$ and $\pdbst$
with the corresponding substacks (over $d$-general  curves). 

The only assertion in \ref{pdg}  that does not hold is the 
existence of a ``Poincar\'e'' line bundle, in (\ref{fine}),
which will never be used in this paper. 
\end{remark}

We are ready to go back to the study of Abel maps.

\begin{prop}
\label{keyrig}
Let $f:\X\to B$ be a regular pencil of $d$-general 
stable curves. Then there exists a canonical map
$$
\ad:\dX^d\la \pf
$$
which restricts to the $d$-th Abel map on the 
generic fiber.
\end{prop}

\begin{proof} We may glue local extensions of 
the $d$-th Abel map of $\X_K$ because they are unique. 
Thus we may assume $f$ is local;
let $X$ be the closed fiber of $f$. In this case, 
the explicit description of $\pf$ is (see \cite{cner}, Cor.~5.14)
\begin{equation}
\label{pfeq}
\pf=\frac{\coprod_{\md \in B^d_X} \Pic^{\md}_f}{\sim_K},
\end{equation}
where, as in \ref{Npre}, ``$\sim_K$" means gluing over the 
generic fiber.

As we know from (\ref{nerdesc}) in \ref{Npre}, 
$\nf$ is described in a very similar way. 
Indeed,
by \cite{cner}, Thm. 6.1, we have a canonical isomorphism
\begin{equation}
\label{picnereq}
\epsilon^d_f:\pf \stackrel{\cong}{\la} \nf.
\end{equation}
To describe it precisely is straightforward: for each 
$\md\in B^d_X$ the restriction of
$\epsilon^d_f$ to $\Pic^{\md}_f$ is the natural isomorphism
$$
\Pic^{\md}_f\stackrel{\cong}\la{} \Pic^{[\md]}_f\subset\nf,
$$
restricting to the identity on the generic fibers. 
The isomorphism $\epsilon^d_f$ is completely
described because, since $X$ is $d$-general, the class map 
$B^d_X\to\dcg^d$ is bijective, by Definition \ref{dgen}.
To conclude, use Lemma \ref{key} and 
(\ref{picnereq}) to define 
$\ad:=(\epsilon^d_f)^{-1}\circ \NN(\alpha^d_K)$.
\end{proof}

We call $\alpha^d_f$ the {\it $d$-th Abel map} of $f$. 
The natural problem now is to 
describe $\alpha^d_f$ as the moduli map
of a balanced line bundle on 
$\pi:\dX^d \times _B \X\to\dX^d$. Since 
$\pf$ is a fine moduli scheme, this should be 
possible. In fact, the proof of 
\cite{caporaso}, Prop. 4.1, p. 621, can be used to 
produce an algorithm 
for determining the necessary twisters we need to tensor 
$\O_{\dX^d\times_B\X}(E^d)$ with to get the balanced line 
bundle.

However the explicit description of this line bundle
 turns out to be difficult to find in general. In Section \ref{open} we will find it for $d=1$.
 In the next 
subsection, \ref{twoc}, 
we will do that for every $d$
in a  special case. 

\begin{nota}\label{twoc}
{\it Two-component curves.}  
Let $X$ be a stable curve with only two irreducible 
components, $C_1$ and $C_2$. Let $g$ be the 
arithmetic genus of $C$ and 
$\delta:=\#C_1\cap C_2$. Set 
(cf. \ref{balrm})
$$
m:=\lceil m_{C_1}(d)\rceil = \Big\lceil \frac{dw_{C_1}}{2g-2}-\frac{\delta}{2}\Big\rceil.
$$
The set $B^d_X$ of balanced multidegrees of $X$ satisfies
\begin{equation}
\label{Bvine}
B^d_X\supseteq\{(m,d-m),(m+1,d-m-1),\dots,
(m+\delta -1,d-m-\delta +1)\},
\end{equation}
with equality if and only if $m_{C_1}(d)$ is not integer, 
if and only if $X$ is $d$-general.

For each integer $a$ such that $0\leq a \leq d$ 
define $r(a)$ to be the  integer determined by
the following two conditions:
$$
0\leq r(a)<\delta\quad\text{ and }\quad 
a-m\equiv r(a) \mod \delta.
$$
Using this notation we have:
\end{nota}

\begin{prop}
\label{vine}Let $X$ be a stable curve with exactly two 
irreducible components, $C_1$ and $C_2$. Let 
$\delta:=\#C_1\cap C_2$. For each regular smoothing
$f:\X\to B$ of $X$, let
$$
\L^{(d)}:=\O_{\dX^d\times_B\X}\Bigl(E^d +
\sum_{a=0}^d\frac{a-m-r(a)}{\delta}
C_1^a\times C_2^{d-a}\times C_1\Bigr)
$$
where, abusing notation, we view 
$$
C_1^a\times C_2^{d-a}\times C_1\subset X^d\times 
X\subset\X^d\times _B \X
$$
as a Cartier divisor
of $\dX^d\times _B \X$, by restriction. Then 
$\L^{(d)}$ is balanced on 
$f_{\dX^d}:\dX^d\times _B\X\to\dX^d$. Furthermore, 
if $X$ is $d$-general, then the $d$-th Abel map 
$\alpha^d_f:\dX^d\to\pf$ 
is the moduli map
of $\L^{(d)}$.
\end{prop}

\begin{proof} Since each 
$C_1^a\times C_2^{d-a}\times C_1$ is supported over 
the closed point of $B$, the line bundle $\L^{(d)}$ 
coincides with $\O_{\dX^d\times_B\X}(E^d)$ over the 
generic point of $B$. This implies that 
$\L^{(d)}$ and $\O_{\dX^d\times_B\X}(E^d)$ 
have the same N\'eron map.

If $X$ is $d$-general, 
$\alpha^d_f$ is defined and coincides 
with the N\'eron map of $\O_{\dX^d\times_B\X}(E^d)$
on $\X_K^d$. On the other hand, if the line bundle 
$\L^{(d)}$ is balanced on 
$\pi:\dX^d\times_B\X\to\dX^d$, there is an associated 
moduli map $\hat\mu_{\L^{(d)}}:\dX^d\to\pf$ by 
\ref{pdg} \eqref{mymap}. Since $\hat\mu_{\L^{(d)}}$ 
coincides with the N\'eron map of $\L^{(d)}$ on $\X_K^d$, 
it follows that $\hat\mu_{\L^{(d)}}=\alpha^d_f$. 
Thus, it suffices to prove that 
$\L^{(d)}$ is balanced.

To verify this, we must compute the 
multidegree of the restriction of
$\L^{(d)}$ to every singular fiber of 
$\pi:\dX ^d\times \X \la \dX^d$. Of course, all singular 
fibers lie over $\dot{X}^d$. So, let $p\in \dot{X}^d$, 
and denote by $X_p$ the fiber $\pi^{-1}(p)$. 
The point $p$ determines a unique pair of 
nonnegative integers $(a_0,b_0)$ such that $a_0+b_0=d$ and 
$p\in C_1^{a_0}\times C_2^{b_0}$.
Then
\begin{equation}
\label{degE}
(E^d\cdot C_1,E^d\cdot C_2) = (a_0, b_0)=(a_0, d-a_0).
\end{equation}
To compute the intersection degrees with $C_1$ and $C_2$ 
of the remaining summands defining $\L^{(d)}$, 
notice first that, 
since $p\in C_1^{a_0}\times C_2^{b_0}$, the only nonzero 
degrees come from the summand indexed by $a=a_0$. 
Now, using 
$$
(C_1^{a_0}\times C_2^{b_0}\times C_1)\cdot (C_1,C_2)=
(-\delta, \delta)
$$
and (\ref{degE}), we get
$$
\mdeg_{X_p} \L^{(d)}=(m+r(a_0),d-m-r(a_0)),
$$
which is balanced because $0\leq r(a_0)<\delta$; see 
(\ref{Bvine}).
\end{proof}

\begin{exstat}
\label{split} Let $X$ be a ``split" curve of arithmetic 
genus $g$, that is, 
$X=C_1\cup C_2$ with $C_i\cong \pr{1}$
and $\#C_1\cap C_2 = g+1$. Then, 
for each $d=1,\dots,g$ and 
any regular smoothing $f:\X\to B$ of $X$,  
the map $\alpha_f^d$ is the moduli map
of $\O_{\dX^d\times_B\X}(E^d)$.

In particular, for any $p_1,\ldots, p_d\in \dot{X}$ 
we have, independently of $f$,
$$
\alpha _f^d(p_1,\ldots, p_d)=\O_X(p_1+\cdots + p_d).
$$
\end{exstat}

A split curve is $d$-general if and only if 
$d\equiv g \mod 2$. Actually, the conclusion of 
Example \ref{split} is  valid 
regardless of $X$ being 
$d$-general, because at any rate 
$\O_{\dX^d\times_B\X}(E^d)$ 
is stably balanced on $f_{\dX^d}$.

\begin{remark}
\label{vinermk} The case of split curves is in a sense 
special. In general, we should expect the 
restriction $\alpha^d_f|_X$ of the 
$d$-th Abel map of Proposition \ref{vine} to depend on the 
choice of smoothing $f$. For a simple 
concrete example of this 
dependence, consider the case $d=2$ and $\delta=2$. Then 
$X$ is stable and $2$-general if $C_1$ and $C_2$ 
have distinct positive arithmetic genera. Suppose $C_1$ 
has smaller genus. Then $m=0$, and thus $r(0)=0$, 
$r(1)=1$, but $r(2)=0$. Let 
$p,\,q\in C_1\smallsetminus C_1\cap C_2$. Then 
$$
\alpha^2_f(p,q)=\O_X(p+q)\otimes\O_\X(C_1)|_X.
$$
Now, as $f$ varies through all smoothings of $X$, 
the restriction 
$\O_\X(C_1)|_X$ varies through all line bundles 
restricting 
to $\O_{C_1}(-r_1-r_2)$ and $\O_{C_2}(r_1+r_2)$, where 
$r_1$ and $r_2$ are the nodes of $X$. 
So $\alpha^2_f(p,q)$ depends on the choice of $f$.

\end{remark}

Denote by $\Sigma^d_g$ the locus in $\mgbar$ of 
curves that are not $d$-general. 
Then $\Sigma^d_g$ is a proper closed 
subset of $\mgbar$; in fact, with the notation of \ref{C94}, $\Sigma^d_g$ is the image via $\phi_d$
of the complement of $U$. As we 
mentioned in \ref{C94}, 
$\Sigma^d_g$ empty if and only if 
$(d-g+1,2g-2)=1$.
Since the rest of our paper is devoted to Abel 
maps for $d=1$,
we conclude this section by describing the locus of curves that are not $1$-general.

\begin{prop}
\label{1gen} Let $g\geq 2$.
\begin{enumerate}[{\rm (i)}]
\item
\label{1gen1}
If $g$ is odd then $\Sigma^1_g$ is empty.
\item
\label{1gen2}
If $g$ is even then $\Sigma^1_g$ is the closure 
in $\mgbar$ of the locus of curves 
$X$ such that $X=C_1\cup C_2$, with $C_1$ and $C_2$ 
smooth of the same genus and $\# C_1\cap C_2$ odd.
\end{enumerate}
\end{prop}

\begin{proof}
If $g$ is odd then
$$
(1-g+1,2g-2)=(g-2,g-1)=1.
$$
So Part \eqref{1gen1} follows; see \ref{dgen}.

For Part \eqref{1gen2}, let $X$ be a stable curve. 
Suppose first that $X$ has the description given 
in \eqref{1gen2}. Then 
      $$\Big(\frac{1-k}{2},\frac{1+k}{2}\Big)\in 
        B^1_X\smallsetminus\tilde B^1_X,$$
where $k:=k_{C_1}=k_{C_2}$. So $X$ is in 
$\Sigma^1_g$. Since $\Sigma^1_g$ 
is closed in $\mgbar$, it contains the closure 
of the locus defined in \eqref{1gen2}.

Suppose now that $X$ is in $\Sigma^1_g$, i.e. 
there is a line bundle $L$ on 
$X$ such that 
$\mdeg L\in B^1_X\smallsetminus\tilde B^1_X$. 
Then there is a connected, 
proper subcurve $Z\subsetneq X$ 
such that either $m_Z(1)$ or $M_Z(1)$ is equal 
to $\deg_ZL$. 
Then both $m_Z(1)$ and $M_Z(1)$ are integers. Thus 
\begin{equation}
\label{B1}
\frac{w_Z}{w}+\frac{k_Z}{2}\in\Z
\end{equation}
where $w:=\deg\omega_X=2g-2$. Now, since 
$X$ is stable,
\begin{equation}
\label{01}
0<\frac{w_Z}{w}<1.
\end{equation}
In particular, $\frac{w_Z}{w}$ is never an integer, 
and thus \eqref{B1} implies that $k_Z$ is odd.
Since $k_Z$ is odd, (\ref{B1}) and (\ref{01}) 
immediately yield $\frac{w_Z}{w}=\frac{1}{2}$.

If $Z'$ is connected, then, since 
$\frac{w_{Z'}}{w}=\frac{1}{2}$ as well, 
we have $g_Z=g_{Z'}$. 
Since both $Z$ and $Z'$ are limits of 
smooth curves, $X$ lies in the closure 
of the locus described in \eqref{1gen2}. 

Thus it remains to show that $Z'$ is also connected. 
Let $Z'_1,\dots,Z'_m$ be the connected 
components of $Z'$. Notice that 
\begin{equation}
\label{sumk}
k_{Z'_1}+\cdots+k_{Z'_m}=k_{Z'}=k_Z.
\end{equation}
Suppose by contradiction that $m>1$. Then 
\begin{equation}
\label{01again}
0<\frac{w_{Z'_i}}{w}<\frac{1}{2}
\end{equation}
for each $i$. Since $\mdeg L\in B^1_X$, we have
      $$\frac{w_{Z'_i}}{w}-\frac{k_{Z'_i}}{2}
        \leq\deg_{Z'_i}L\leq
	\frac{w_{Z'_i}}{w}+\frac{k_{Z'_i}}{2}$$
for each $i$. Using \eqref{01again}, we get
      $$\frac{1-k_{Z'_i}}{2}\leq\deg_{Z'_i}L\leq
        \frac{k_{Z'_i}}{2}$$
for each $i$. Summing up, and using \eqref{sumk}, 
we get
      $$\frac{m-k_Z}{2}\leq\deg_{Z'}L\leq
        \frac{k_Z}{2}.$$
Now, since $\deg L=1$, we must have 
\begin{equation}
\label{bdzL}
\frac{2-k_Z}{2}\leq\deg_ZL\leq\frac{2-m+k_Z}{2}.
\end{equation}

Suppose first that 
$\deg_ZL=M_Z(1)$, that is, $\deg_ZL=(1+k_Z)/2$. 
Then \eqref{bdzL} implies $1+k_Z\leq 2-m+k_Z$, and 
hence $m\leq 1$, a contradiction.

Finally, suppose $\deg_ZL=m_Z(1)$, that is, 
$\deg_ZL=(1-k_Z)/2$. Then \eqref{bdzL} implies 
$2-k_Z\leq 1-k_Z$, a contradiction as well.
\end{proof}

\section{Modular interpretation of the first Abel map}
\label{open}
The following diagram represents the families we 
shall deal with in this section,
starting from a regular pencil of stable curves
$f:\X\to B$:
\begin{equation}
\begin{array}{lcccr}
\dX\times_B\X & \ha&\X^2_B&\la &\X \\
\  \  \downarrow{\dpi} &&\downarrow{\pi} &&\downarrow{f}\\
\   \dX & \ha&\X&\la &B \\
\end{array}
\end{equation}
where $\pi$ is the projection onto the first factor.

We denote by $\Delta \subset \X^2_B$ the diagonal. 
Its restriction to $\dX\times_B\X$ is a Cartier divisor. 
Denote 
by $\O_{\dX\times_B\X}(\Delta )$ the associated line 
bundle. We may view $\O_{\dX\times_B\X}(\Delta )$ 
as a family of degree-1 line bundles on the fibers 
of $\dpi$. Recall that the first Abel map of the 
generic fiber of $f$ is the moduli map of the 
restriction of $\O_{\dX\times_B\X}(\Delta )$; see 
\ref{abelmap}. We want to interpret the 
first Abel map $\alpha^1_f$, 
defined in Proposition \ref{keyrig}, 
as the moduli map of a 
balanced line bundle on $\dpi$, which will necessarily 
be a (possibly trivial) twist of 
$\O_{\dX\times_B\X}(\Delta )$.

In fact, we shall see that 
$\O_{\dX\times_B\X}(\Delta )$ fails to be balanced 
over points of a singular fiber $X$ of $f$ only if 
$X$ has a separating node. To fix this, 
we will tensor $\O_{\dX\times_B\X}(\Delta )$ by 
twisters supported on so-called ``tails''.

We need a few preliminary results which hold for any curve 
$X$, possibly having  singularities other than nodes. 
For the sake of future applications of the 
techniques developed in this paper, 
from now until \ref{ttt}, and in \ref{Qtwister}, 
\ref{septree} and \ref{preAXinj}, 
we shall be in this more general situation,
i.e. $X$ will be any (reduced, connected and projective) 
curve over an algebraically closed field.

Let $r$ be a node of $X$ and 
$X^{\nu}_r\to X$ be the normalization of $X$ at $r$ only. 
  If  $X^{\nu}_r$ is not connected, $r$ is called  a {\it separating node} of $X$.

\begin{defi}\label{deft}
Let $X$ be a curve of arithmetic genus $g$. 
A proper subcurve 
$Q\subsetneq X$ will be called a {\it tail} of $X$ if 
$Q$ intersects the complementary subcurve $Q'$
in a separating node $r$ of $X$. We say that $Q$ is 
{\it attached to $r$} or that $r$ {\it generates $Q$}. 
A tail $Q$ of $X$ 
will be called {\it small} if $g_Q<g/2$ and 
{\it large} if $g_Q>g/2$. Let
$$
\QQ(X):=\{ Q\subset X : Q\text{ is a small tail of } X \}.
$$
If $X$ has no separating node, for instance if $X$ 
is smooth, then $\QQ(X)=\emptyset$.
\end{defi}

If $r$ is a separating node  of $X$, 
then $X^{\nu}_r$ has   two 
connected components, isomorphic to the two tails generated by $r$; hence
  every 
tail  is connected. 

For every tail $Q\subset X$ we
have that $g=g_Q+g_{Q'}$. So, at least one of the two 
tails attached to a separating
node has arithmetic genus at most $g/2$.  
If the curve $X$ is stable and $1$-general, it follows 
from Proposition~\ref{1gen} that 
no tail of $X$ can have genus 
equal to $g/2$, or in other words that 
every tail of $X$ is either small or large.

\begin{remark}\label{intail}
Let $r$ be a separating node of $X$ generating the tails
$Q$ and $Q'$.
If  $Z\subset X$ is a connected
subcurve not containing $r$, then 
$Z$ is entirely contained in either $Q$ or $Q'$.
\end{remark}

\begin{lemma}\label{ZZ}
Let $X$ be a curve and $Q_1$ and $Q_2$ two tails of $X$. 
Then
$$
Q_1\cup Q_2=X
\quad\text{or}\quad 
Q_1\cap Q_2=\emptyset
\quad\text{or}\quad
Q_1\subseteq Q_2
\quad\text{or}\quad
Q_2\subseteq Q_1.
$$
\end{lemma}

\begin{proof} 
For each $i=1,2$ let $r_i$ be the separating node 
of $X$ generating $Q_i$. If $r_1=r_2$ then 
either $Q_1=Q_2'$, and hence $Q_1 \cup Q_2 = X$, 
or $Q_1=Q_2$. So we may assume that $r_1\neq r_2$.

Thus $r_1\not\in Q_2$ or $r_1\not\in Q_2'$. Suppose first 
that $r_1\not\in Q_2$. Since $Q_2$ is connected, 
either $Q_2\subset Q_1$ or $Q_2\subset Q_1'$ by 
Remark \ref{intail}. If $Q_2\subset Q'_1$ then 
$Q_1\cap Q_2\subseteq Q_1\cap Q_1'=\{r_1\}$, and 
hence $Q_1\cap Q_2=\emptyset$ because 
$r_1\not\in Q_2$. So 
either $Q_2\subset Q_1$ or $Q_1\cap Q_2=\emptyset$.

The case where $r_1\not\in Q_2'$ is treated similarly. 
In this case, 
either $Q_2'\subset Q'_1$, and hence $Q_1\subset Q_2$, 
or $Q_1'\cap Q_2'=\emptyset$, and hence $Q_1\cup Q_2=X$.
\end{proof}

\begin{lemma}\label{ctype}
Let $X$ be a curve, and $Q$ a tail of $X$. 
Then, for any two line bundles $L_1$ on $Q$ and $L_2$ 
on $Q'$, 
there is, up to isomorphism, a unique line bundle $L$ 
on $X$ such that $L|_Q\cong L_1$ and 
$L|_{Q'}\cong L_2$.
\end{lemma}

\begin{proof} Let $r$ be the separating node of $X$ 
to which $Q$ is attached. For each isomorphism 
$\mu:L_1|_{\{r\}}\to L_2|_{\{r\}}$, let $L$ be 
the kernel of the composition
$$
\phi_\mu:L_1\oplus L_2 \la 
L_1|_{\{r\}}\oplus L_2|_{\{r\}} 
\stackrel{\tilde\mu}{\longrightarrow} 
L_2|_{\{r\}},
$$
where $\tilde\mu:=(-\mu,1)$. Since $\phi_\mu$ is 
surjective, $L\neq L_1\oplus L_2$. Also,
since $\mu$ is an isomorphism, $L\neq L_1(-r)\oplus L_2$ 
and $L\neq L_1\oplus L_2(-r)$. Since $r$ is a node 
of $X$, it follows that $L$ is a line bundle, and 
$L|_Q\cong L_1$ and 
$L|_{Q'}\cong L_2$. 

Conversely, if $N$ is a line bundle on $X$ for which 
there are isomorphisms $\lambda_1:N|_Q\to L_1$ and 
$\lambda_2:N|_{Q'}\to L_2$, then $N$ is the 
kernel of $\phi_\mu$, where
$\mu:=\lambda_2|_{\{r\}}\circ\lambda^{-1}_1|_{\{r\}}$.
 
Finally, if $\mu':L_1|_{\{r\}}\to L_2|_{\{r\}}$ 
is another isomorphism, the kernel of $\phi_\mu$ 
is carried isomorphically to the kernel of 
$\phi_{\mu'}$ by the autoomorphism
       $$(a,1):L_1\oplus L_2\la L_1\oplus L_2,$$
where $a$ is the unique scalar such that $\mu=a\mu'$.
\end{proof}

\begin{nota}\label{ttt} {\it Twisters on tails.} 
Let $X$ be a curve. By Lemma~\ref{ctype}, 
for each tail $Q$ of $X$ 
there is a unique, up to isomorphism, 
line bundle on $X$ whose 
restrictions to $Q$ and $Q'$ are $\O_Q(-r)$ and 
$\O_{Q'}(r)$, where $r$ is the separating node 
of $X$ generating $Q$. Denote this bundle by 
$\O_X(Q)$. 

For each formal sum $\sum a_QQ$ of tails $Q$ with 
coefficients $a_Q\in\Z$, set
$$
\O_X({\textstyle\sum}a_QQ):=\bigotimes
\O_X(Q)^{\otimes a_Q}.
$$ 
If $X$ is a nodal curve, and a closed fiber of 
a regular pencil $f:\X\to B$, then
\begin{equation}
\label{OXX}
\O_{\X}(\textstyle\sum a_QQ)|_X\cong\O_X(\sum a_QQ).
\end{equation}
So {\it twisters supported on tails do not 
depend on the chosen regular pencil}. To check 
\eqref{OXX} it 
is enough to observe that, for each tail $Q$ of $X$, 
since 
$$
\O_{\X}(Q)|_Q\cong\O_Q(-r)\text{ and }
\O_{\X}(Q)|_{Q'}\cong\O_Q(r),
$$
Lemma~\ref{ctype} implies that 
$\O_\X(Q)|_X\cong\O_X(Q)$.
\end{nota}

To state the main result of this section we need 
some notation, similar to the one used in 
Proposition \ref{vine}. 
Let $f:\X\to B$ be a regular pencil. 
Let $Z$ be a subcurve of $X$, where $X\subset\X$ is 
a singular fiber of $f$. Then $Z$ is a divisor 
of $\X$. Now, the restriction $\pi_{Z}$ of 
the first projection $\pi:\X^2_B\to\X$ over $Z$ 
is the trivial family
$$
\X^2_B\supset Z\times X \stackrel{\pi_{Z}}{\la }Z.
$$
Thus, for any other subcurve $Z_1\subseteq X$, the 
product $Z\times Z_1$ can be viewed as a Weil 
divisor of $\X^2_B$. Now, since the open subscheme 
$\dX\times_B\X \subset\X^2_B$ is regular, 
the restriction of $Z\times Z_1$ to it is a Cartier 
divisor. Let 
$\O_{\dX\times_B\X}(Z\times Z_1)$ denote the 
associated line bundle. Using this notation, we 
have an explicit description of the map 
$
\alpha_f^1:\dX \la \pfu 
$
defined in Proposition \ref{keyrig}.

\begin{thm}
\label{AN1} Let $f:\X\to B$ be a regular pencil of 
stable curves. Then the line bundle
$$
\L^{(1)}:=\O_{\dX\times_B\X}\Bigl(\Delta +
\sum_{\stackrel{b\in B}{Q \in \QQ( X_b)}}Q\times Q 
\Bigr).
$$
is balanced on $\dpi:\dX\times_B\X\to\dX$. Furthermore, 
assume that the fibers of $f$ are $1$-general. Then 
the following two statements hold.
\begin{enumerate}[{\rm (i)}]
\item\label{mod} The morphism $\alpha_f^1:\dX\to P^1_f$ 
is the moduli map of $\L^{(1)}$. 
\item\label{uni} If ${\mathcal M}$ is a balanced line 
bundle on $\dpi:\dX\times_B\X\to\dX$
having $\au$ as moduli map, then 
${\mathcal M}\cong\L^{(1)}$, up to pullbacks from $\dX$.
\end{enumerate} 
\end{thm}

\begin{remark}
But for Part \eqref{uni}, the theorem
generalizes to curves that are not 1-general. See 
\ref{nogen} and Remark~\ref{Q1fix}.
\end{remark}

\begin{proof} The divisor $\sum Q\times Q$ is 
entirely supported on a union of closed fibers of 
$\X\times_B\X\to B$. By Lemma \ref{key}, the 
N\'eron maps of $\L^{(1)}$ and $\O(\Delta)$ are equal. 
Now, if the fibers of $f$ are 1-general, $\alpha^1_f$ 
agrees with the N\'eron map of $\O(\Delta)$ 
on $\X_K$. On the other hand, if 
$\L^{(1)}$ is balanced on $\dpi:\dX\times_B\X\to\dX$, 
there is an associated moduli map 
$\hat\mu_{\L^{(1)}}:\dX\to P^1_f$ by \ref{pdg} \eqref{mymap}. 
Since $\hat\mu_{\L^{(1)}}$ 
coincides with the N\'eron map of $\L^{(1)}$ on $\X_K$, 
it follows that $\hat\mu_{\L^{(1)}}=\alpha^1_f$. Thus, 
to prove Part \eqref{mod} it is enough to prove that 
$\L^{(1)}$ is balanced on $\dpi$.

To prove Part (\ref{uni}), it is also enough to show 
that $\L^{(1)}$ is balanced, since $P^1_f$ is a fine 
moduli scheme; see \ref{pdg} \eqref{fine}.

Let us now prove that $\L^{(1)}$ is indeed balanced. 
We need only check this on each singular fiber of $f$, 
whence we may assume $f$ is local. Let 
$X$ be the closed fiber. It suffices to  consider 
the singular fibers of the first projection
$$
\dpi:\dX\times_B\X \la \dX,
$$
which are all isomorphic to $X$. 
Let $p\in X$ be a nonsingular point, and 
set $L_p^{(1)}:=\L^{(1)}|_{\dpi^{-1}(p)}$. Then 
\begin{equation}
\label{reseq}
L_p^{(1)}\cong \O_X(p)\otimes
\O_X(\sum_{\stackrel{Q\in \QQ(X)}{p\in Q}}Q).
\end{equation}
We conclude by  Lemma~\ref{balp} (\ref{balpt}),
observing that, since $X$ is stable, 
a semibalanced line bundle on $X$
is necessarily balanced; see \ref{balrm} (\ref{balrm1}).
\end{proof}

The next two lemmas are needed to finish the proof of 
Theorem \ref{AN1}.

\begin{lemma}\label{Qtwister}
Let $X$ be a curve, and $Z$ a connected,
proper subcurve. Let
\begin{equation}\label{chain}
Q_1\subset Q_2\subset\cdots\subset Q_{n-1}\subset Q_n
\end{equation}
be a chain of tails of $X$,
and let $r_i$ be the
separating node of $X$ generating $Q_i$ for each 
$i=1,\dots,n$. Then
$$-1\leq\deg_Z\O_X(\textstyle\sum Q_i)\leq 1.$$
Furthermore, the extremes are attained if and only if there 
is a unique $j$ such that $r_j\in Z\cap Z'$. In this 
case, the lower bound is attained if $Z\subseteq Q_j$, 
and the upper bound is attained if $Z\subseteq Q'_j$.
\end{lemma}

\begin{proof}
If $r_\ell\not\in Z$ for any $\ell=1,\dots,n$, then
$\deg_Z\O_X(\textstyle\sum Q_i)=0$.
Suppose now that $Z$ contains at least
one $r_\ell$. Let $i$ and $j$ be the smallest
and greatest integers
such that $r_i\in Z$ and $r_j\in Z$, respectively.
Since $Z$ is connected,
$Z$ contains as well all the irreducible
components of $X$ containing
$r_{i+1},\dots,r_{j-1}$. In particular, 
$r_\ell\not\in Z\cap Z'$ for any $\ell=i+1,\dots,j-1$. 

If $Z\cap Z'$ contains both $r_i$ and $r_j$ or neither 
of them, $\deg_Z\O_X(\textstyle\sum Q_i)=0$. If 
$Z\cap Z'$ contains $r_i$ but not $r_j$, then 
$\deg_Z\O_X(\textstyle\sum Q_i)=1$ and 
$Z\subseteq Q'_i$. At last, if $Z\cap Z'$ contains $r_j$ 
but not $r_i$, then 
$\deg_Z\O_X(\textstyle\sum Q_i)=-1$ and 
$Z\subseteq Q_j$. 
\end{proof}

\begin{lemma}
\label{balp} 
Let $X$ be a semistable curve, and $p$ a nonsingular 
point. Then the following two statements hold.
\begin{enumerate}[{\rm (i)}]
\item
\label{balpp}The line bundle $\O_X(p)$ is 
semibalanced if and only if $p$ does not belong to 
any small tail of $X$.
\item
\label{balpt}
The line bundle
$$L_p^{(1)}:=\O_X(p)\otimes
\O_X(\sum_{\stackrel{Q\in \QQ(X)}{p\in Q}}Q)$$
is semibalanced.
\end{enumerate}
\end{lemma}

\begin{proof} The ``if part'' of \eqref{balpp} is a 
consequence of \eqref{balpt}, as the sum of tails in 
\eqref{balpt} is zero when $p$ does not belong to any 
small tail. As for the ``only-if part'', 
recall from \ref{notcurve} that
\begin{equation}\label{adjw}
w_Z=2g_Z-2+k_Z
\end{equation}
for every connected proper subcurve 
$Z\subset X$. In particular, 
\begin{equation}\label{smt<g-1}
w_Q<g-1\text{ for every small tail $Q$ of $X$.}
\end{equation} 
So, if $p$ is contained in a small tail $Q$, 
then 
$$
M_Q(1)<1=\deg_Q\O_X(p).
$$
Hence the 
Basic Inequality \eqref{BIeq} is 
not satisfied for $Q$. So $\O_X(p)$ is not semibalanced. 

We need only prove \eqref{balpt} now. First, 
since $X$ 
is semistable, $w_Z\geq 0$ for every subcurve 
$Z\subseteq X$; see \ref{notcurve}. 
As a consequence, 
\begin{equation}\label{w>0}
w_{Z_1}\leq w_{Z_2}
\text{ for all subcurves $Z_1$ and $Z_2$ of $X$ with 
$Z_1\subseteq Z_2$.}
\end{equation}

Let $Q_1,\dots,Q_n$ be the small tails of $X$ 
containing $p$, and $r_1,\dots,r_n$ their 
generating  nodes. Since $w_{Q_i}+w_{Q_j}<2g-2$ by 
\eqref{smt<g-1}, 
we have $Q_i\cup Q_j\neq X$ for each $i$ and $j$. 
By Lemma \ref{ZZ}, up to reordering, we 
may assume that
$$
Q_1\subset Q_2\subset\cdots\subset Q_{n-1}\subset 
Q_n.
$$

Let $N:=\O_X(\sum_1^n Q_i)$; so 
$L_p^{(1)}=\O_X(p)\otimes N$. 
Let $Z$ be any connected, proper subcurve of 
$X$. Then $\deg_Z N\geq -1$ 
by Lemma \ref{Qtwister}, and hence 
$\deg_ZL_p^{(1)}\geq -1$. As pointed out in 
Remark \ref{balrm} (\ref{balrm3}), we need only show 
that $\deg_ZL_p^{(1)}\geq m_Z(1)$.

First, suppose 
$\deg_ZN=-1$. By Lemma \ref{Qtwister}, there is $j$ 
such that $r_j\in Z$ and $Z\subseteq Q_j$. 
By~\eqref{smt<g-1} 
and \eqref{w>0},
\begin{equation}\label{wz<g-1}
w_Z\leq w_{Q_j}<g-1,
\end{equation}
and hence $m_Z(1)\leq 0$. Thus, if $p\in Z$, 
     $$\deg_Z L_p^{(1)}=0\geq m_Z(1).$$

Suppose $p\not\in Z$. Then $Z\neq Q_j$. Since 
$r_j\in Z$, either $k_Z\geq 3$ or $Z$ is a tail of 
$Q_j$. Now, if $Z$ were a tail of $Q_j$, then 
$\overline{Q_j\smallsetminus Z}$ would be a tail of 
$X$ contained in $Q_j$, whence a small 
tail. Since $p\in\overline{Q_j\smallsetminus Z}$, we have 
$\overline{Q_j\smallsetminus Z}=Q_i$ for some $i<j$, 
or $Z=\overline{Q_j\smallsetminus Q_i}$. But then 
$\deg_ZN=0$, a contradiction. 
Thus $k_Z\geq 3$. In particular, $Z$ is not an 
exceptional component of $X$. It follows now from 
\eqref{wz<g-1} that $m_Z(1)<-1$, and hence
      $$\deg_Z L_p^{(1)}\geq -1>m_Z(1).$$

Second, suppose $\deg_ZN\geq 0$. 
Then $\deg_ZL_p^{(1)}\geq 0$.
By \eqref{w>0} we have
that $w_Z\leq w_X=2g-2$. So, if $Z$ is not a large 
tail of $X$, then $m_Z(1)\leq 0$, and hence
      $$\deg_Z L_p^{(1)}\geq 0\geq m_Z(1).$$
On the other hand, suppose that 
$Z$ is a large tail. At any rate, $m_Z(1)\leq 1/2$. 
Thus, if $p\in Z$,
\begin{equation}\label{dZ>1}
\deg_Z L_p^{(1)}\geq 1\geq 1/2\geq m_Z(1).
\end{equation}
Finally, suppose $p\not\in Z$. Then $p$ lies on $Z'$, 
which is a small 
tail of $X$. Thus $Z'=Q_j$ for some $j$, and 
hence $Z=Q'_j$. It follows that 
$\deg_ZN=1$, and hence 
\eqref{dZ>1} holds as well.
\end{proof}

Let $X$ be a 1-general stable curve. 
Let $\dot{X}:=X\smallsetminus X_{\text{sing}}$. 
For any regular smoothing $f$ of $X$, let 
$$
\alpha^1_X:=\alpha^1_f|_{\dot{X}}:\dot{X}\la P_X^1.
$$ 
The notation is not ambiguous
by the following consequence of Theorem~\ref{AN1}.

\begin{cor}
\label{AX}
Let $X$ be a 1-general stable curve. Then 
$\alpha^1_X$ does not depend on $f$. In fact, 
for each nonsingular point $p\in X$ we have
$$
\alpha^1_X(p)=\O_X(p)\otimes
\O_X(\sum_{\stackrel{Q\in\QQ(X)}{p\in Q}}Q).
$$
\end{cor}

\begin{proof} 
The expression of $\alpha^1_X(p)$ follows from (\ref{reseq}) in the proof of Theorem~\ref{AN1}.
By \ref{ttt} the map $\alpha^1_X$ does not depend on $f$.
\end{proof}

{\it If $X$ is free from separating nodes
then $\alpha^1_X$ is injective.} 
This follows immediately from Lemma \ref{preAXinj} 
below. The same lemma will be used in the proof 
of Proposition~\ref{A1cinj}, a more general and 
precise statement. For the lemma and the 
proposition, the definition below is used.

\begin{defi}\label{septree} 
Let $X$ be a curve. A rational, smooth 
component $C$ of $X$ is called a {\it separating line} 
if $C$ intersects $\overline{X\smallsetminus C}$ 
in separating nodes of $X$. More generally, a 
connected subcurve $Z\subseteq X$ of arithmetic genus 0 
is called a {\it separating tree of lines} if 
$Z$ intersects $\overline{X\smallsetminus Z}$ in 
separating nodes of $X$.
\end{defi}

\begin{nota}\label{septreermk} 
Let $X$ be a curve, and $Z\subsetneq X$ a proper 
connected subcurve such that $Z$ intersects 
$\overline{X\smallsetminus Z}$ in separating nodes 
of $X$. Then the connected components of 
$\overline{X\smallsetminus Z}$ are tails of $X$. 
In addition, if $r$ is a separating node of $Z$, 
then $r$ is a separating node of $X$. 

A curve of arithmetic genus 0 is a 
curve of compact type, i.e. a nodal curve with every node 
separating, whose irreducible components are 
smooth and rational. So, if $Z$ is a separating tree of 
lines, every node of $Z$ is a separating node of $Z$, 
and hence of $X$. It follows that every 
connected subcurve of $Z$ is also a separating tree of 
lines. 
In particular, every irreducible component of $Z$ is a separating line.
\end{nota}

We shall later need the following lemma.

\begin{lemma}
\label{preAXinj} Let $X$ be a curve, and $p$ and $q$ 
distinct nonsingular points of $X$. Let $C\subseteq X$ 
be the irreducible component containing $p$. 
Then there is an isomorphism $\O_X(p)\cong\O_X(q)$ 
if and only if $C$ contains $q$ and is a separating line 
of $X$.
\end{lemma}

\begin{proof} Assume first that $C$ contains $q$ and 
is a separating line 
of $X$. Since $C$ is smooth and 
rational, $\O_C(p)\cong\O_C(q)$. We may thus assume 
$C\neq X$. Since $C$ meets 
$C':=\overline{X\smallsetminus C}$ in separating 
nodes, applying Lemma \ref{ctype} a few times, we can 
show that a line bundle 
on $X$ is uniquely determined by its restrictions to 
$C$ and to $C'$. Since $\O_X(p)$ and $\O_X(q)$ restrict 
to isomorphic line bundles on $C$ and to the trivial 
line bundle on $C'$, it follows that 
$\O_X(p)\cong\O_X(q)$.

Conversely, suppose $\O_X(p)\cong\O_X(q)$. 
Since $\O_X(p)$ has degree 1 on $C$, so 
has $\O_X(q)$, and hence $q\in C$ as well. Now, since 
$\O_X(p)\cong\O_X(q)$, in particular
$\O_C(p)\cong\O_C(q)$. Since $p\neq q$, it follows 
from \cite{AK}, Thm. 8.8, p.~108 that $C\cong\pr{1}$. 

If $C=X$ we are done. Suppose thus that $C\neq X$, and 
let $C':=\overline{X\smallsetminus C}$. Also, suppose 
by contradiction that $C\cap C'$ is 
not made of separating nodes of $X$. 
Then there is a connected subcurve $Z\subseteq C'$ 
such that $C\cap Z$ is a scheme of length at 
least 2. 

Since $X$ is connected, the restriction 
$\tau:H^0(X,\O_X(p))\to H^0(C,\O_C(p))$ is injective. 
But it is not surjective. Indeed, if a nonconstant
$\sigma\in H^0(C,\O_C(p))$ could be 
extended to $\tilde{\sigma}\in H^0(X,\O_X(p))$, 
then $\tilde{\sigma}$ would have to be constant on $Z$ and 
hence $\sigma$ would be constant on $C\cap Z$. 
Since $C\cong \pr{1}$, this is impossible, $\sigma$ being a nonconstant 
section of $\O_{\pr{1}}(1)$
 So $\tau$ is injective, but not 
surjective, and hence $h^0(X,\O_X(p))=1$. Since 
$\O_X(p)\cong\O_X(q)$, it follows that $p=q$, 
an absurd.
\end{proof}

\section{Completing the first Abel map}
\label{close}

The main result of this section is 
Theorem \ref{A1c}. We shall prove it first in a 
simpler case in Proposition \ref{A1nosep}, 
where we have a neater statement concerning the 
modularity, see \ref{A1mod}. 

As in Section \ref{open}, 
certain basic results of this section hold 
in more generality for curves having singularities 
other than nodes. Apart from the notation set in 
\ref{blow} below, these results are concentrated 
in \ref{ressing}.

\begin{nota}\label{blow}
Let $X$ be a  curve. 
For each node $r\in X$, let $X^{\nu}_r\to X$ denote 
the partial normalization of $X$ at $r$, and let 
$\lad_r$ be  the
curve obtained by adding to $X^{\nu}_r$
a smooth rational curve $E_r$ 
connecting the two branches over $r$. Thus
$$
\lad_r= X^{\nu}_r\cup E_r,
$$
and there is a natural surjection
$\sigma_r:\lad_r\to X$ such that 
$\sigma_r(E_r)=\{r\}$, and such that 
$\sigma_r$ is an isomorphism away from $E_r$.

Assume now that $X$ is a 1-general stable curve, and 
let $r$ be a nonseparating node of $X$. Let 
$\pren\in E_r\subset \lad_r$ such that $\pren$ is a 
nonsingular point of $\lad_r$. Then the line bundle 
$
\O_{\lad_r}(\pren)\in \Pic^1\lad_r
$
is balanced by Lemma \ref{balp}, and hence 
determines a point of $\pub\smallsetminus \pu$; see 
\ref{pdg} \eqref{mymap}. This point 
does not depend 
on the choice of $\pren$; see \ref{equiv}. Thus 
we shall denote it by $\LN$.
\end{nota}

\begin{prop}
\label{A1nosep} Let $f:\X \to B$ be a regular pencil of 
1-general stable curves free from separating nodes. 
Then $\alpha^1_f:\dX\to P^1_f$ extends to an 
injection
$$
\aub:\X \la \pfub
$$
such that $\aub(r)=\LN\in\pub$ for each node $r$ of 
each closed fiber $X$ of $f$.
\end{prop}

\begin{remark}
\label{A1mod}
More precisely, the proof will show that $\aub$ is 
the moduli map of the line bundle 
$\O_\Y(\tilde{\Delta})$, where $\Y\to\X^2_B$ is 
a partial resolution of singularities,
and $\tilde{\Delta}$ is the proper transform in 
$\Y$ of the diagonal $\Delta$; see \ref{ressing}.
\end{remark}

\begin{proof}
Denote by $\rho:\Y\to\X^2_B$ the partial resolution 
of singularities described in \ref{ressing}, from where 
we take some of the properties mentioned below. The 
map $\rho$ is an isomorphism away from the 
points $(r,r)$ for $r\in\X\smallsetminus\dX$. On the 
other hand, if $r\in\X\smallsetminus\dX$, then 
$\rho^{-1}(r,r)$ is a copy of $\pr{1}$.
In addition, composing $\rho$ with the first 
projection $\pi$,
$$
\Y\stackrel{\rho}{\la} \X^2_B\stackrel{\pi}{\la}\X,
$$
we obtain a family of quasistable curves 
$\Y\to\X$ having $\pi:\X^2_B\to\X$ as stable model. 

For each closed fiber $X$ of $f$, and each 
$r\in X_{\text{sing}}\subset \X$, let $Y_r$ 
be the fiber of $\pi\circ\rho$ over $r$. Then 
$Y_r=\lad_r$, where $\lad_r$ is as defined in 
\ref{blow}. On the other hand, 
each fiber of $\pi\circ\rho$ 
over $\dX$ is the same as the corresponding fiber 
of $\pi$.

Let $\tilde{\Delta}\subset \Y$ be the proper 
transform of $\Delta$. By Property \ref{ressing} 
\eqref{rsdiag}, the map 
$\rho$ restricts to an isomorphism between 
$\tilde{\Delta}$ and ${\Delta}$. Also, 
$\tilde{\Delta}$ meets each fiber $Y_r=\lad_r$ of 
$\pi\circ\rho$ over 
$\X\smallsetminus\dX$ 
transversally at a nonsingular point $\pren$ contained
in the exceptional component $E_r$. 

To prove that $\alpha^1_f$ extends, we prove two claims: 
first, that 
$\O_\Y(\tilde{\Delta})$ is 
balanced on $\pi\circ\rho:\Y\to\X$, so 
it induces a morphism
$\aub:\X\to\pfub$, its moduli map; and second, 
to show that $\aub$ extends $\au$, that the 
restriction of $\O_\Y(\tilde{\Delta})$ to
the each fiber of $\pi\circ\rho$ 
over $\dX$ is isomorphic to the 
corresponding restriction of $\L^{(1)}$, 
whose moduli map is $\au$ by Theorem \ref{AN1}. 

We may now assume that $f$ is local. Let $X$ be its closed 
fiber. 
For each $r\in X_{\text{sing}}$, since $\tilde{\Delta}$ 
intersects $Y_r$ transversally at $\pren$, we have
\begin{equation}\label{OYD}
\O_\Y(\tilde{\Delta})|_{Y_r}\cong\O_{\lad_r}(\pren),
\end{equation}
which is balanced by Lemma \ref{balp}. 
In addititon, for each nonsingular point $p\in X$ we have
$$
\O_\Y(\tilde{\Delta})|_{Y_p}\cong\O_{X}(p)
$$
which is balanced and isomorphic to the 
corresponding restriction of $\L^{(1)}$, also 
by Lemma \ref{balp}. Therefore $\O_\Y(\tilde{\Delta})$ 
induces a moduli map
$$
\aub:\X \la \overline{P^1_f}
$$
which extends $\au$. Notice that \eqref{OYD} also shows 
that $\aub(r)=\LN$ for each $r\in X_{\text{sing}}$. 

To show that $\aub$ is 
injective it suffices to consider singular points 
of $X$, by Corollary \ref{AX} and by the fact that
$\aub(r)\in\pub\smallsetminus\pu$ for every node 
$r\in X$. Now, if $r\in X_{\text{sing}}$, then $\aub(r)$ 
represents a balanced line bundle on $\lad_r$. Hence, 
two different nodes $r$ and $r'$ of $X$ are mapped to two 
points of $\pub$ corresponding to balanced line bundles
on different quasistable curves, namely $\lad _r$ and
$\lad_{r'}$. Thus $\aub(r)\neq \aub (r')$; see 
\ref{equiv}.
\end{proof}

\begin{nota}\label{ressing}
{\it Resolution of singularities.}
In the proof of Proposition \ref{A1nosep} we used a 
partial resolution of singularities of $\X^2_B$ 
which we are
now going to describe in detail, and in more 
generality. 

Let $f:\X\to B$ be a regular pencil. 
The threefold $\X^2_B$ is singular at the points 
$(r_1,r_2)$, where $r_1$ and $r_2$ are (not necessarily 
distinct) singular points of the same closed fiber of $f$. 

Let $X$ be a closed fiber of $f$, and $r_1$ and $r_2$ 
nodes of $X$. Since $f$ is regular, 
locally around $r_i$ the 
surface $\X$ is formally equivalent to the 
surface in $\A^3$ given by the equation
$x_iy_i=t$, where $t$ denotes a local parameter of 
$B$ at the closed point covered by $X$. 
Pulling back these 
local equations to $\X^2_B$ under the two projection maps 
$\X^2_B\to\X$, and abusing of the same notation, 
we get that $\X^2_B$ is 
formally equivalent, locally around $(r_1,r_2)$, 
to the threefold in $\A^5$ 
given the equations
\begin{displaymath}
\left\{ \begin{array}{l}
x_1y_1=t,\\
x_2y_2=t.\\
\end{array}\right.
\end{displaymath}
If $r_1=r_2$, then the diagonal 
$\Delta\subset\X^2_B$ contains $(r_1,r_2)$, 
and we may assume that it is given locally 
around $(r_1,r_2)$ by
\begin{displaymath}
\left\{ \begin{array}{l}
x_1y_1=t,\\
x_2=x_1,\\
y_2=y_1.\\
\end{array}\right.
\end{displaymath}

Locally around $(r_1,r_2)$ we may eliminate $t$, 
and view $\X^2_B$ as the cone $C\subset\A^4$ 
over the smooth quadric 
in $\pr{3}$ given by $x_1y_1=x_2y_2$. Also, 
if $r_1=r_2$, we may view $\Delta$ as the 
plane $D\subset\A^4$ given by $x_2=x_1$ and 
$y_2=y_1$. Notice that $C$ is singular only at the 
origin. To resolve this 
singularity we need only blow up a plane in $C$ containing 
the origin. Any plane will do, but let us blow up the 
plane 
 given 
by $x_1=x_2=0$. The blowup is 
the nonsingular 
threefold $\tilde C\subset\pr{1}\times \A^4$ given 
by the equations
\begin{displaymath}
\left\{ \begin{array}{l}
\xi_2x_1=\xi_1x_2,\\
\xi_1y_1=\xi_2y_2,\\
\end{array}\right.
\end{displaymath}
where $\xi_1,\xi_2$ are homogeneous coordinates 
of $\pr{1}$. The blowup $\gamma:\tilde C\to C$ 
is isomorphic to $C$ away from the origin. 
In addition, the fiber 
$F$ over the origin is given by 
$x_1=x_2=y_1=y_2=0$, and 
hence is isomorphic to $\pr{1}$. 

The exceptional divisor $E$ of the blow up 
$\tilde C$ is given by $x_2=0$ where $\xi_2\neq 0$, 
and $x_1=0$ where $\xi_1\neq 0$. In particular, 
$F\subset E$. Now, since $\xi_2x_1=\xi_1x_2$, 
summing the divisor given by $\xi_1=0$ to $E$ 
we get the principal divisor given by $x_1=0$. 
Thus $E\cdot F=-1$.

Suppose $r_1=r_2$. Then $\gamma^{-1}(D)$ is given by 
$x_1(\xi_1-\xi_2)=y_2(\xi_1-\xi_2)=0$ where 
$\xi_1\neq 0$, and by 
$x_2(\xi_1-\xi_2)=y_1(\xi_1-\xi_2)=0$ where 
$\xi_2\neq 0$. Thus $\gamma^{-1}(D)$ is the 
union of the Cartier divisor given by $\xi_1=\xi_2$ 
and the fiber $F$. The strict transform 
$\tilde D$ of $D$ is thus a Cartier divisor 
intersecting $F$ transversally at a point. 

For $i=1,2$, let $\phi_i:\tilde C\to\A^2$ be the 
composition of $\gamma$ with the projection onto 
the plane with coordinates $x_i,y_i$. 
Its fiber over the origin is given 
by $x_1=y_1=\xi_1x_2=\xi_2y_2=0$. 
It is the union of $F$ and the affine lines 
$N_1$, given by $x_1=y_1=\xi_1=y_2=0$, and $N_2$, given by
$x_1=y_1=\xi_2=x_2=0$. The lines $N_1$ and $N_2$ do not 
meet, and $F$ intersects each $N_i$ transversally 
at a single point. Also, $\phi_2$ maps $N_1$ and 
$N_2$ isomorphically onto the lines $y_2=0$ and 
$x_2=0$, respectively.

The exceptional divisor $E$ contains $N_2$, and 
intersects $N_1$ transversally. Since $\xi_1\neq 0$ 
on $N_2$, we have $E\cdot N_2=0$. If $r_1=r_2$, 
the strict transform $\tilde D$ does not meet 
either $N_1$ or $N_2$, and intersects $F$ 
transversally. 

(An analogous description holds if we 
reverse the roles of $\phi_1$ and $\phi_2$.)

We will now consider the global picture. Let 
${\mathcal I}_\Delta$ denote the ideal sheaf of the 
diagonal $\Delta\subset\X^2_B$, and 
let $\check{\mathcal I}_\Delta$ denote 
the dual sheaf, i.e.
      $$\check{\mathcal I}_\Delta:=
        Hom({\mathcal I}_\Delta,\O_{\X^2_B}).$$
Since ${\mathcal I}_\Delta$ is a sheaf of ideals, 
$\check{\mathcal I}_\Delta$ is a sheaf of fractional 
ideals of $\X^2_B$. A piece of notation: for each open 
subscheme 
$U\subseteq\X^2_B$ and each sheaf of fractional ideals 
$\mathcal M$ of $U$, consider its powers 
${\mathcal M}^n$, and form the 
associated sheaf of Rees algebras:
      $${\mathcal R}({\mathcal M}):=\O_U\oplus{\mathcal M}
        \oplus{\mathcal M}^2\oplus\cdots\oplus
	{\mathcal M}^n\oplus\cdots.$$
Set $\Y:=\text{Proj}_{\X^2_B}({\mathcal R}
(\check{\mathcal I}_\Delta))$, and let 
$\rho:\Y\to\X^2_B$ be the structure map. 

We may view $\rho$ as a blowup. In fact, for any open 
subscheme $U\subseteq\X^2_B$ over which there is an 
embedding 
$\iota:\check{\mathcal I}_\Delta|_U\to{\mathcal L}$ 
into an invertible sheaf ${\mathcal L}$, we may view 
$\rho:\rho^{-1}(U)\to U$ as the blowup of $U$ along the 
closed subscheme $V\subseteq U$ whose sheaf of ideals 
${\mathcal I}_{V|U}$ satisfies 
$\iota(\check{\mathcal I}_\Delta|_U)=
{\mathcal I}_{V|U}{\mathcal L}$. In other words, $\iota$ 
induces an isomorphism over $U$:
      $$\text{Proj}_{U}({\mathcal R}
        ({\mathcal I}_{V|U}))\longrightarrow
        \text{Proj}_{U}({\mathcal R}
        (\check{\mathcal I}_\Delta)|_U).$$
In the same vein, for each invertible sheaf of ideals 
${\mathcal J}\subseteq\O_U$ we have that 
$Hom({\mathcal I}_\Delta|_U,{\mathcal J})=
\check{\mathcal I}_\Delta|_U{\mathcal J}$, and hence 
we obtain a canonical isomorphism over $U$:
      $$\text{Proj}_{U}({\mathcal R}
        (\check{\mathcal I}_\Delta)|_U)\longrightarrow
        \text{Proj}_{U}({\mathcal R}
	(Hom({\mathcal I}_\Delta|_U,{\mathcal J}))).$$

Since ${\mathcal I}_\Delta$ is invertible 
away from the points $(r,r)$ for 
$r\in\X\smallsetminus\dX$,
it 
follows from the above description that 
$\rho$ is an isomorphism away from these same points. 
In addition, 
around the points $(r,r)$, where $r$ 
is a node of a closed fiber of $f$, 
the map $\rho$ is formally 
equivalent to the blowup described above, because
      $$\text{Hom}_A((x_1-x_2,y_1-y_2),(x_1-x_2))
        =(x_1,x_2),\text{ where }
	A:=\frac{{\mathbb C}[[x_1,x_2,y_1,y_2]]}
	{(x_1y_1-x_2y_2)}.$$
Then all of the properties above, verified locally, 
yield global properties of $\rho$. Indeed, assume 
that the fibers of $f$ are nodal. (It would actually 
be enough to assume that the fibers are Gorenstein.) 
Then, recalling that $\pi:\X^2_B\to\X$ denotes the
first projection, the following statements hold:
\begin{enumerate}[{\bf (A)}]
\item
\label{rsgen} The composition 
$$\pi\circ\rho:\Y \la \X$$
is a family of curves whose fiber 
$Y_r$ over a point $r$ of a closed fiber 
$X$ of $f$ is $X$, 
if $r$ is nonsingular, 
and $\lad_r$, described in {\rm\ref{blow}}, 
if $r$ is a node.
\item
\label{rsdiag}
 Let $\tilde{\Delta}\subset\Y$ denote the 
proper transform of $\Delta$. 
For each node $r$ of each closed fiber $X$ of $f$, the 
transform $\tilde{\Delta}$ intersects the fiber 
$Y_r$ transversally at a point lying in the 
exceptional component $E_r=\rho^{-1}(r,r)$.
\item
\label{rsint}
 Let $Q$ be a tail of a closed fiber $X$ of $f$, 
and $r$ the node of $X$ 
generating $Q$. Let
$$
\widetilde{Q^2}:=\rho^{-1}(Q\times Q).
$$ 
Then $\widetilde{Q^2}$ is a Cartier 
divisor of $\Y$ containing $E_r$. Furthermore,
$$
\widetilde{Q^2}\cdot E_r=-1,\quad 
\widetilde{Q^2}\cdot\hat Q=-1,\quad\text{and}\quad 
\widetilde{Q^2}\cdot\hat Q'=1,
$$
where, using the notation in \ref{blow}, 
$\hat Q:=\sigma_r^{-1}(Q)$
and $\hat Q':=\overline{\lad_r\smallsetminus \hat Q}$, 
i.e. $\hat Q$ is the tail of $\lad_r$ mapping to $Q$ and 
containing $E_r$, and $\hat Q'$ is the complementary 
tail.
\end{enumerate}
\end{nota}

We may now generalize Proposition \ref{A1nosep}.

\begin{thm}
\label{A1c}
Let $f:\X\to B$ be a regular pencil of 1-general 
stable curves. Then there exists a morphism
$$
\aub:\X \la \overline{P^1_f}
$$
extending $\au:\dX\to P^1_f$. 
If $r$ is a node of a closed fiber $X$ of $f$, then 
$\aub(r)\in\pu$ if and only if $r$ is a separating node 
of $X$. 
\end{thm}

\begin{remark}
The result extends to curves that are not 1-general. See \ref{nogen} and \ref{Q1fix}.
\end{remark}

\begin{proof} As in the proof of Proposition \ref{keyrig} 
we may work locally around each singular fiber. So, assume 
$f$ is local, and let $X$ denote its closed fiber.

The new difficulty with respect to 
Proposition \ref{A1nosep}
is that, if $X$ has separating nodes, $\au$ is the 
moduli map of a nontrivial ``twist'' of the diagonal by 
Theorem \ref{AN1}, and thus the 
same must hold for its completion. Fortunately, 
however, the  divisors we need for the ``twist'' 
are already present in the partial resolution 
of singularities $\rho:\Y\to\X^2_B$ described 
in \ref{ressing}.

Namely, let $Q_1,\dots,Q_m$ be all the small tails 
of $X$. Let $\tilde\Delta\subset \Y$ be the 
strict transform 
of $\Delta$, and set 
$\widetilde{Q_i^2}:=\rho^{-1}(Q_i\times Q_i)$ 
for $i=1,\dots,m$. As seen in \ref{ressing}, 
all the $\widetilde{Q_i^2}$  and $\tilde\Delta$ are 
Cartier divisors. Define the line bundle
\begin{equation}
\label{M}
{\mathcal M}:=\O_\Y\Bigl(\tilde\Delta+ 
\widetilde{Q_1^2}+\cdots+\widetilde{Q_m^2}\Bigr).
\end{equation}
We claim that ${\mathcal M}$ is semibalanced 
on the composition $\pi\circ\rho:\Y\to\X$ of $\rho$ 
with the first projection $\pi:\X^2_B\to\X$. Once 
the claim is proved, we may let 
$\aub:\X \la \overline{P^1_f}$ be the moduli map of 
${\mathcal M}$; see \ref{pdg} \eqref{mymap}. 

To prove the claim, first observe that 
$\rho$ is an isomorphism over $\dX\times_B\X$, 
whence 
$$
{\mathcal M}|_{\dX\times_B\X}\cong
\O_{\dX\times_B\X}(\Delta+Q_1\times Q_1+\cdots+
Q_m\times Q_m),
$$
which is balanced, by Theorem \ref{AN1}, and defines 
$\au$. Thus, 
once $\mathcal M$ is shown to be semibalanced, 
we have that $\aub|_{\dX}=\au$.

Now, let $r\in X_{\text{sing}}$. The fiber 
$Y_r:=(\pi\circ\rho)^{-1}(r)$ is 
equal to $\lad_r$ by 
Property~\ref{ressing}~(\ref{rsgen}). Also, 
$\tilde\Delta$ intersects $Y_r$ transversally at a 
point $\pren$ of the exceptional component 
$E_r=\rho^{-1}(r,r)$, by Property \ref{ressing} 
(\ref{rsdiag}).

For each $i=1,\dots, m$, let 
$r_i$ be the separating 
node of $X$ generating $Q_i$. Let 
$$
\hat Q_i:=\sigma _r^{-1}(Q_i)\subset\lad_r,
$$ 
and let $\hat Q'_i$ be its complement in $\lad_r$. 
Then $\hat Q_i$ is a small tail 
of $\lad_r$ dominating $Q_i$, and containing $E_r$ if and 
only if $r\in Q_i$. If $r=r_i$ then also 
$\overline{\hat Q_i\smallsetminus E_r}$ 
is a small tail of $\lad_r$. These 
are all the 
small tails of $\lad_r$: the subcurves 
$\hat Q_1,\,\dots,\,\hat Q_m$, together with 
$\overline{\hat Q_i\smallsetminus E_r}$ in case $r=r_i$.

For each $i=1,\dots,m$, the subscheme 
$Q_i\times Q_i\subset\X^2_B$ is a Cartier divisor 
away from $(r_i,r_i)$. Identifying $Y_r$ with $\lad_r$, 
we claim that 
      \begin{equation}\label{OYQ}
	\O_\Y(\widetilde{Q_i^2})|_{Y_r}\cong
	\begin{cases}
	  \O_{\lad_r}&\text{if $r\not\in Q_i$,}\\
	  \O_{\lad_r}(\hat Q_i)&
	  \text{if $r\in Q_i$.}
	\end{cases}
       \end{equation}
In fact, if $r\not\in Q_i$, then $\widetilde{Q_i^2}$ 
does not meet $Y_r$, and hence \eqref{OYQ} holds. 
Suppose now that $r\in Q_i$. Recall that 
$\hat{Q_i}$ is a tail of $\lad_r$. Let $s_i$ denote its 
generating node. If $r\neq r_i$ then, since 
$Q_i\times Q_i$ is a Cartier divisor of $\X^2_B$ at 
$(r,r_i)$, we have
    \begin{equation}\label{OYQs}
      \O_\Y(\widetilde{Q_i^2})|_{\hat Q_i}\cong
      \O_{\hat Q_i}(-s_i)
\quad\text{and}\quad
      \O_\Y(\widetilde{Q_i^2})|_{\hat Q'_i}\cong
      \O_{\hat Q'_i}(s_i).
    \end{equation}
The same restrictions are achieved with 
$\O_{\lad_r}(\hat Q_i)$. Thus \eqref{OYQ} follows 
from Lemma \ref{ctype}. Finally, if $r=r_i$ then 
\eqref{OYQs} still holds, 
by Property \ref{ressing} \eqref{rsint}, and hence 
\eqref{OYQ} follows in the same way. The proof of 
\eqref{OYQ} is complete.

Now, notice that $Q_i$ contains $r$ if and only if 
$\hat Q_i$ contains $\pren$. In addition, if $r=r_i$ 
then $\pren\not\in\overline{\hat Q_i\smallsetminus E_r}$. 
Since 
$\hat Q_1,\,\dots,\,\hat Q_m$, and 
$\overline{\hat Q_i\smallsetminus E_r}$ if $r=r_i$, 
are all the small tails of $\lad_r$, 
it follows from \eqref{OYQ} that
$$
{\mathcal M}|_{Y_r}\cong\O_{\lad_r}(\pren)
\otimes\O_{\lad_r}
(\sum_{\stackrel{Q\in\QQ(\lad_r)}{\pren\in Q}}Q),
$$
which is semibalanced by Lemma \ref{balp}. Our claim 
is proved, and thus we finish the proof of 
the existence of $\aub$.

To prove the second statement of the theorem, 
it suffices to prove that for any node
$r\in X$ we have
\begin{displaymath}
\deg _{E_r}{\mathcal M} =\left\{ \begin{array}{ll}
1 &\text{ if $r$ is not separating,}\\
0 &\text{ otherwise.}\\
\end{array}\right.
\end{displaymath}

To prove this, notice that, if $r\neq r_i$ then 
$\widetilde{Q_i^2}\cdot E_r=0$, whereas if 
$r=r_i$ then $\widetilde{Q_i^2}\cdot E_r=-1$ by 
Property \ref{ressing} \eqref{rsint}. 
Since at any rate 
$\tilde\Delta\cdot E_r=1$, the degree of 
${\mathcal M}|_{E_r}$ is 1, unless 
$r=r_i$ for some $i$, in which case the degree is 0. 
Since $X$ is 1-general, each separating node of $X$ 
generates a small tail, and hence is equal to $r_i$ for 
some $i$. So $\aub(r)\in\pu$ if and only if $r$ is a 
separating node.
\end{proof}

\begin{example}\label{exct}
Let $X$ be a curve of compact type with two components, 
$C_1$ and $C_2$. Then $C_1$ and $C_2$ are smooth, 
and $C_1\cap C_2=\{r\}$, where $r$ is the unique node 
of $X$. Assume $g_{C_1}<g_{C_2}$. 
Then $\lad_r = C_1\cup E\cup C_2$ and
$\QQ(\lad_r)=\{ C_1, C_1\cup E\}$, 
where $E=\pr{1}$. The line bundle 
${\mathcal M}$ in the proof of
Theorem \ref{A1c}, whose moduli map is $\aub$, satisfies
$$
{\mathcal M}= \O_{\Y}\Bigl( {\tilde{\Delta}}+ 
\widetilde{C_1^2} 
\Bigr).
$$
In this case, it is easy to describe the completed 
Abel map. First notice that there is a 
canonical isomorphism 
$\pub \cong \Pic^0C_1\times \Pic ^1C_2$, essentially 
by Lemma \ref{ctype}. Hence, a point 
$\ell \in \pub$ 
is represented by a pair $(L_1,L_2)$ with 
$L_i\in \Pic C_i$. For $i=1,2$ let $q_i\in C_i$ lying 
above $r$. Then 
\begin{displaymath}
\aub (p)=\left\{ \begin{array}{ll}
(\O_{C_1}(p-q_1), \O_{C_2}(q_2)) &\text{ if } p\in C_1, \\
(\O_{C_1}, \O_{C_2}(p)) &\text{ if } p\in C_2.\\
\end{array}\right.
\end{displaymath}
In particular, $\aub(r)=(\O_{C_1}, \O_{C_2}(q_2))$. 
Thus, composing $\aub|_{C_1}$ with the projection
$\Pic^0C_1\times \Pic ^1C_2 \to \Pic^0C_1$
we obtain the classical Abel--Jacobi map of $C_1$ 
with base point $q_1$, i.e. 
$$
\begin{array}{ccc}
C_1&\la &\Pic^0C_1\\
p&\mapsto&\O_{C_1}(p-q_1).\\  
\end{array}
$$
The analogous composition for $C_2$ gives the first 
Abel map $C_2\to \Pic^1C_2$.
\end{example}

Let $X$ be a 1-general stable curve, and $f$ a 
regular smoothing of $X$. The restriction $\aub|_{\dot X}$ 
coincides with $\alpha^1_f|_{\dot X}$, whence 
does not depend on $f$ by Corollary~\ref{AX}. So 
$\aub|_X$ does not depend on $f$ either, and we may set
$\aXub:=\aub|_X$.

\begin{nota} {\it Explicit description of the complete Abel map.} Let $X$ be a 1-general stable curve 
and $p\in X$. We shall now explicitely describe 
$\aXub(p)$, in formulas (\ref{psm}) and (\ref{psi}) below,
following the proof
of  Theorem~\ref{A1c}.

First some notation.
Let $P_1,\dots,P_m$ be all the small 
tails of $X$ containing $p$. 
(The unusual naming of the tails 
using ``$P$" rather than ``$Q$"
is to match
the notation of the proof of 
Proposition~\ref{A1cinj}.) 
By Lemma \ref{ZZ} we can
write 
$$
P_m\subset P_{m-1}\subset\dots\subset P_2\subset P_1.
$$
Set
$Z_i:=\overline{P_i-P_{i+1}}$ for each 
$i=1,\ldots, m-1$ and $Z_m:=P_m$,
so that $P_1=\cup_1^mZ_i$. Also, put $Q:=P_m'$, 
a large tail of $X$. Hence
$$
X=P_1\cup Q=\cup_1^mZ_i\cup Q.
$$
Let $r_1,\dots,r_m$ be the separating nodes of 
$X$ generating $P_1,\dots,P_m$. Notice that 
$Z_i\cap Z_i' =\{r_i, r_{i+1}\}$ if 
$i=1,\ldots, m-1$,
and $Z_m\cap Z_m' =\{r_m\}$. 
Therefore each of the $Z_i$ and $Q$ meets the 
complementary curve in separating
nodes of $X$. Hence, by iterated use of 
Lemma~\ref{ctype}, to give a line bundle on $X$ 
it suffices to give its
restrictions to all the $Z_i$ and to $Q$.

We are now ready to describe $\aXub(p)$
if $p$ is a nonsingular point or
a separating node of $X$ (in which case of 
course $p=r_m$). Recall that by
Theorem~\ref{A1c}, $\aXub(p)$ corresponds to a line bundle on $X$. We have
\begin{equation}
\label{psm}
\aXub(p)=\{\O_{Q}(r_1), \O_{Z_1}(r_2-r_1),\ldots, \O_{Z_{m-1}}(r_m-r_{m-1}),
\O_{Z_m}(p-r_m)\}.
\end{equation}

Now, suppose that $p$ is a nonseparating node of $X$. Then we know that $\aXub(p)$
corresponds to a line bundle on $\lad _p$. 
Let $E\subset \lad _p$ be the exceptional component of $\lad_p$, and let $\widetilde{Z_m}$ denote 
the normalization of $Z_m$ at $p$ only.
Keeping the above notation we have
$$
\lad _p=Q\cup Z_1\cup\cdots\cup Z_{m-1}\cup \widetilde{Z_m}\cup E.
$$
Now, recall from \ref{equiv} that $\aXub(p)$ is 
uniquely determined by a line bundle 
$L$, of degree $0$, on the complementary curve of 
$E$; that is, arguing as above, by the string of 
the restrictions of $L$ to $Q,\,Z_1,\,
\ldots,\, Z_{m-1},\, \widetilde{Z_m}$.
We have
\begin{equation}
\label{psi}
\aXub(p)=\{\O_{Q}(r_1), \O_{Z_1}(r_2-r_1),\ldots,\O_{Z_{m-1}}(r_m-r_{m-1}),
\O_{\widetilde{Z_m}}(-r_m)\}
\end{equation}
\end{nota}

\

\begin{prop}\label{A1cinj} 
Let $X$ be a 1-general stable 
curve. Let $p$ and $q$ be distinct points of $X$. Then 
$\aXub(p)=\aXub(q)$ if and only if $p$ and $q$ belong to 
the same separating tree of lines of $X$.
\end{prop}

A similar result for the Abel--Jacobi map to the (degree-0) Jacobian 
is proved by B.Edixhoven in \cite{edix}, Prop. 9.5.
His statement (necessarily) excludes the case where 
$p$ or $q$ is a nonseparating node, since there the target space of
the map is a noncompactified N\'eron model.

\begin{proof} Suppose first that $p$ and $q$ belong to a 
separating tree of lines of $X$, call it $Z$. Since $Z$ is 
connected, to prove that $\aXub(p)=\aXub(q)$ it is 
enough to consider the case where 
$p$ and $q$ are nonsingular points 
of $X$ in the same irreducible component $C$ of $Z$. 
Now, $C$ is a separating line of $X$; 
see \ref{septreermk}. Thus $\O_X(p)\cong\O_X(q)$ by 
Lemma \ref{preAXinj}. Since $p$ and $q$ lie on the 
same component, it follows that 
$\alpha^1_X(p)=\alpha^1_X(q)$, and hence 
$\aXub(p)=\aXub(q)$.
 
Conversely, suppose $\aXub(p)=\aXub(q)$. 
We claim that $p$ and $q$ are either nonsingular points or 
separating nodes of $X$. Indeed, suppose by contradiction, 
and without loss of generality, that $p$ is a 
nonseparating node of $X$. Then 
$\aXub(p)\in\pub\smallsetminus\pu$ by Theorem \ref{A1c}. 
Since $\aXub(p)=\aXub(q)$, it follows from 
Theorem \ref{A1c} as well that $q$ is also 
a nonseparating node of $X$. 
However, $\aXub(p)$ and $\aXub(q)$ 
correspond to balanced line bundles on different 
quasistable curves, $\lad_p$ and $\lad_q$. 
So $\aXub(p)\neq\aXub(q)$; see \ref{equiv}. The 
contradiction proves the claim. 

Since $p$ and $q$ are nonsingular or 
separating nodes of $X$, both $\aXub(p)$ and 
$\aXub(q)$ are line bundles on $X$. 
Let $P_1,\dots,P_m$ be the small 
tails containing $p$ and $Q_1,\dots,Q_n$ the small tails 
containing $q$. (We may have $m=0$ or 
$n=0$.) It follows from Lemma \ref{ZZ}, as in the 
proof of Theorem \ref{AN1}, that, up to reordering the 
tails,
$$
P_m\subset P_{m-1}\subset\dots\subset P_2\subset P_1
\quad\text{and}\quad
Q_n\subset Q_{n-1}\subset\dots\subset Q_2\subset Q_1.
$$
Set $P_0:=Q_0:=X$. Let 
$r_1,\dots,r_m$ be the separating nodes of $X$ generating 
$P_1,\dots,P_m$, and $s_1,\dots,s_n$ those generating 
$Q_1,\dots,Q_n$. In addition, 
set $P_{m+1}:=Q_{n+1}:=\emptyset$, and 
put $r_{m+1}:=p$ and $s_{n+1}:=q$. 

We may assume $m\leq n$, without loss of generality. 
Let $i$ be the largest nonnegative integer such that 
$i\leq m$ and $P_j=Q_j$ for $j=0,1,\dots,i$. Then 
also $r_j=s_j$ for $j=0,1,\dots,i$. We claim that 
$P_{i+1}\cap Q_{i+1}=\emptyset$. Indeed, if $i=m$ then 
$P_{i+1}$ is already empty. Suppose $i<m$. If 
$P_{i+1}\subseteq Q_{i+1}$, then 
$$
P_{i+1}\subseteq Q_{i+1}\subset Q_i=P_i,
$$ 
and hence $P_{i+1}=Q_{i+1}$, contradicting the 
maximality of $i$. In a similar way, 
$Q_{i+1}\not\subseteq P_{i+1}$. Since 
$P_{i+1}\cup Q_{i+1}\neq X$, because $P_{i+1}$ and 
$Q_{i+1}$ are small tails, it follows from 
Lemma \ref{ZZ} that $P_{i+1}\cap Q_{i+1}=\emptyset$, 
proving our claim. In particular, $r_{i+1}\neq s_{i+1}$.

As $P_i=Q_i$, we may consider 
$Y:=\overline{P_i\smallsetminus(P_{i+1}\cup Q_{i+1})}$. 
As $P_{i+1}$ and $Q_{i+1}$ do not meet, their 
union cannot be $P_i$, a connected subcurve of $X$. 
Thus $Y$ is a subcurve of $X$. It is also connected, 
being either equal to, or a tail of, 
$\overline{P_i\smallsetminus P_{i+1}}$, which 
in turn
is either equal to, or a tail of, $P_i$, a tail of $X$.

Since $Y\subseteq\overline{P_i\smallsetminus P_{i+1}}$, 
the restriction of $\aXub(p)$ to $Y$ is 
$\O_Y(r_{i+1}-r_i)$. Analogously, 
$\aXub(q)$ restricts to $\O_Y(s_{i+1}-s_i)$. 
Since $\aXub(p)=\aXub(q)$ and $r_i=s_i$, it 
follows that $\O_Y(r_{i+1})\cong\O_Y(s_{i+1})$. Since 
$r_{i+1}\neq s_{i+1}$, by Lemma~\ref{preAXinj} applied 
to the curve $Y$, we see 
that $r_{i+1}$ and $s_{i+1}$ are contained in a 
separating line $C$ of $Y$. Since 
$Y\cap Y'$ is made of separating nodes of $X$, so 
is $C\cap C'$; see \ref{septreermk}. 
In other words, $C$ is a separating line of $X$.

For $\ell=1,\dots,m-i$ let 
$Y_\ell:=\overline{P_{i+\ell}-P_{i+\ell+1}}$. Then 
$\aXub(p)$ restricts to 
$\O_{Y_\ell}(r_{i+\ell+1}-r_{i+\ell})$. 
On the other hand, since 
$Y_\ell\subset\overline{Q_i\smallsetminus Q_{i+1}}$ for 
each $\ell$, but neither $s_i\in Y_\ell$ nor 
$s_{i+1}\in Y_\ell$, we have that $\aXub(q)$ restricts 
to the trivial bundle $\O_{Y_\ell}$. Applying 
Lemma \ref{preAXinj} to the curve $Y_\ell$, we get 
that $r_{i+\ell}$ and 
$r_{i+\ell+1}$ are contained in a separating line 
$C^p_\ell$ of $Y_\ell$. 
As before, $C^p_\ell$ is also a separating line of $X$. 

Similarly, for each $\ell=1,\dots,n-i$ the points 
$s_{i+\ell}$ and $s_{i+\ell+1}$ are contained in a 
separating line $C^q_\ell$ of $X$. The union of all 
the separating lines, namely,
$$
C^p_{m-i},\,\dots,\,C^p_2,\,C^p_1,\,C,\,
C^q_1,\,C^q_2,\,\dots,\,C^q_{n-i},
$$
is a separating tree of lines containing $p$ and $q$. 
\end{proof}

\begin{nota}
\label{nogen} {\it Curves that are not $1$-general.}
We conclude by discussing the case of stable curves 
which are not $1$-general. Recall that such curves
 form a proper closed subset of $\mgbar$, nonempty 
if and only if $g$ is even, and their 
combinatorial structure is described in 
Proposition~\ref{1gen}.

What kind of complications occur for curves that are not $1$-general, or not $d$-general?

The stack $\pdbst$ introduced in \ref{pdg} in the 
case $(d-g+1,2g-2)=1$ is constructed as a quotient
stack, i.e. $\pdbst:=[H_d/G]$; notation as in \ref{C94}.
The same definition can be given for each $d$, 
obtaining in this way a quotient stack
$[H_d/G]$. However, when non-$d$-general curves appear, 
this stack presents some pathologies.

More precisely, recall from \ref{C94} 
that the scheme-theoretic quotient $H_d/G$
is endowed with a natural surjective morphism 
$\phi_d: H_d/G\to\mgbar$. 
The open subset of $\mgbar$ over which the quotient 
map $\pi_d:H_d \to H_d/G$ is a geometric quotient
is exactly the locus of $d$-general curves. 
The problem is that, as soon as  $\pi_d:H_d\to H_d/G$ 
fails to be a geometric quotient, the following 
pathologies occur:
\begin{enumerate}[(i)]
\item $[H_d/G]$ fails to be a Deligne--Mumford stack. 
\item The natural map of stacks $[H_d/G]\to \mgbst$ 
fails to be representable.
\item N\'eron models are not 
parametrized
by $[H_d/G]$.
\end{enumerate}

However, when studying Abel maps, 
we can still obtain some results.
Since the stack $\pdbst$ 
behaves badly, let us consider 
the scheme $\pdgbar:=H_d/G$ introduced in \ref{C94}.
As we mentioned above, there is always a surjective 
morphism $\phi_d:\pdgbar \to \mgbar$. 
By \cite{caporaso}, Thm. 6.1, p. 641, $\pdgbar$ is 
an integral projective scheme. It is also normal, 
being a GIT-quotient of $H_d$, which is
nonsingular by \cite{caporaso}, Lemma 2.2, p. 609. 

Although 
$\pdgbar$ 
is not a coarse moduli space, 
not even away from curves with nontrivial 
automorphisms, $\pdgbar$ 
does satisfy useful functorial
properties.

Thus, let $f:\X \to B$ be a regular pencil of 
stable curves. Let $B\to\mgbar$ be the associated map, 
and define
$$
\pfb := B\times_{\mgbar}\pdgbar.
$$
If $X$ is a closed fiber of $f$, denote by $\pXb$ the 
corresponding fiber of $\pfb$ over $B$. 
As mentioned above, $\pfb$ may fail to contain the 
N\'eron model $N_f^d$. However, a functorial property 
holds: the moduli property given in 
\ref{pdg} (\ref{mymap}) holds exactly as stated.
More precisely, to any semibalanced line bundle $\L$
on a family of semistable curves $\Y \to T$ 
having $\X_T\to T$ as stable model, where $T$ is any 
$B$-scheme, we can associate a canonical moduli map
$\bal_{\L}:T\la \pfb$; see \cite{caporaso}, 
Prop. 8.1, p. 653.

The main weakness, when nongeneral curves are present, 
is that different balanced line bundles on the same 
quasistable, or even stable, curve
may be mapped to the same point in $\pfb$.
Let us give an example of this behavior with regard to 
Abel maps.
\end{nota}

\begin{example}
Let $X=C_1\cup C_2$ be a curve of compact type 
as in Example~\ref{exct}. However, assume now 
that $C_1$ and $C_2$ have the same genus. Then $X$ 
is not $1$-general. As before, let $r$ be the unique 
node of $X$, and let $q_i$ be the point of $C_i$ lying 
over $r$ for $i=1,2$. 

We shall now exhibit three 
nonequivalent balanced line bundles that correspond to 
the same point of $\pub$. Notice that, by Lemma 
\ref{ctype}, to give a line bundle on a curve of compact 
type is equivalent to give a line bundle on each 
irreducible component of the curve.

Let $p\in C_1\smallsetminus\{q_1\}$. 
Our first line bundle is 
$L_1\in \Pic X$ corresponding to the pair
$$
(\O_{C_1}(p), \O_{C_2}).
$$
Let $Y:=\lad _r$; so $Y=C_1\cup E\cup C_2$, where $E$ 
is the exceptional component. Our second line bundle is 
$L_2\in \Pic Y$ corresponding to the triple 
$$
(\O_{C_1}(p-q_1), \O_E(q_E), \O_{C_2}),
$$
where $q_E$ is the point of $E$ glued to $q_1$ in $Y$. 
Finally, our third line bundle is $L_3\in \Pic X$ 
corresponding to the pair
$$
(\O_{C_1}(p-q_1), \O_{C_2}(q_2)).
$$

We leave out all proofs, referring to 
\cite{caporaso}, 7.2, Example 2, p. 645 for more details.
\end{example}

The above example shows that, if $f:\X\to B$ is a regular 
pencil with a non-$1$-general fiber $X$, then 
$\pfub$ and $\pub$ are not coarse 
moduli schemes for balanced line bundles.

However, we can still get a map 
$\aub: \X \to \pfub$ restricting to the classical 
Abel map of $\X_K$, by using our modular interpretation 
of the Abel map. In fact, essentially the same line 
bundle $\mathcal M$ given in \eqref{M} can be used to 
produce a moduli map $\aub: \X \to \pfub$. Most of the 
results in Sections \ref{open} and \ref{close} hold,
provided we change one definition, as explained in Remark 
\ref{Qfix} below.

\begin{nota}
\label{Qfix} 
{\it Small tails, again.}
Let $X$ be a stable curve of arithmetic genus $g$. 
Taking into account the case where $X$ is not 1-general, 
we need to adjust the definition of the set $\QQ(X)$ in 
\ref{deft}.
Suppose that $X$ has a separating node that 
generates two tails 
$Q$ and $Q'$ of equal genus. 
It is easy to see that, if such a node exists, 
then it is unique. (A $1$-general curve will never 
admit such a node by Proposition \ref{1gen}.)
We must add to the set $\QQ(X)$ of small tails of $X$ 
either $Q$ or $Q'$, thus making an arbitrary choice 
between $Q$ and $Q'$, which nonetheless turns out to be 
completely irrelevant. 

So, {\it $\QQ(X)$ is defined as the set of all small tails 
of $X$ together with  one tail of genus $g/2$,
if any such tail exists}.
\end{nota}

\begin{remark}
\label{Q1fix}
The 
following results of the paper hold with 
essentially the 
same proof, as long as we use the modified definition of 
$\QQ(X)$ for stable curves $X$ of \ref{Qfix}:
\begin{enumerate}[(i)]
\item Theorem~\ref{AN1}, excluding Part \eqref{uni}.
\item Corollary~\ref{AX}.
\item Theorem~\ref{A1c}.
\end{enumerate}
What will certainly fail is the possibility to 
interpret the Abel map in a unique way.
In other words, if $f:\X\to B$ is a regular pencil, 
an extension $\aub:\X\to\pfub$ of the Abel map of $\X_K$ 
is obtained as 
the moduli map of a semibalanced line bundle,
however the line bundle is not uniquely determined.
\end{remark}

\

\noindent L.\ Caporaso   (caporaso@mat.uniroma3.it)\\
Dipartimento di Matematica, Universit\`a Roma Tre\\
Largo S.\ L.\ Murialdo 1 \  \  00146 Roma - Italy\\

\noindent E.\ Esteves (esteves@impa.br)\\
IMPA\\
Est. D. Castorina 110 \  \  22460-320 
Rio de Janeiro - Brazil
\end{document}